\documentclass[fleqn,final]{elsart_hal}
\usepackage[english]{babel}

\usepackage{times,epsf,amsmath,amstext,amsfonts}

\newcommand{\bfF}{\text{${\bf F}$}}
\newcommand{\bfE}{\text{${\bf E}$}}

\newcommand{\bfx}{\text{${\bf x}$}}

\newcommand{\bfpsi}{\text{\boldmath $\psi$}}

\newcommand{\bfa}{\text{${\bf a}$}}

\newcommand{\bfD}{\text{${\bf D}$}}
\newcommand{\bfA}{\text{${\bf A}$}}
\newcommand{\bfu}{\text{${\bf u}$}}

\newcommand{\bfB}{\text{${\bf B}$}}
\newcommand{\bfC}{\text{${\bf C}$}}

\newcommand{\bfP}{\text{${\bf P}$}}

\newcommand{\bfdelta}{\text{\boldmath $\delta$}}

\newcommand{\bfb}{\text{${\bf b}$}}
\newcommand{\bfc}{\text{${\bf c}$}}

\newcommand{\bfV}{\text{${\bf V}$}}

\newcommand{\bfzero}{\text{\boldmath ${\bf 0}$}}

\newcommand{\rd}{\text{${\rm d}$}}
\newcommand{\sd}{\text{${\, \rm d}$}}
\newcommand{\rD}{\text{${\rm D}$}}

\begin{document}

\begin{frontmatter}

\title{Numerical simulation of spray coalescence in an Eulerian framework: direct quadrature method of moments and multi-fluid method}

\author{R.\ O.\ Fox\corauthref{vis}}
\address{Department of Chemical and Biological Engineering, 2114 Sweeney Hall,
Iowa State University, Ames, IA 50011-2230, U.S.A. - rofox@iastate.edu}
\corauth[vis]{The present research was conducted during two visits of R.O. Fox in France,
supported by Ecole Centrale Paris (Feb.--July 2005 and June 2006).}

\author{F.\ Laurent\corauthref{cor}}
\corauth[cor]{Corresponding author.} 
and
\author{M.\ Massot}
\address{Laboratoire d'{\'E}nerg{\'e}tique Mol{\'e}culaire et Macroscopique, Combustion, UPR CNRS 288 - {\'E}cole Centrale Paris, Grande Voie des Vignes, 92295 Ch{\^a}tenay-Malabry, France - frederique.laurent@em2c.ecp.fr, marc.massot@em2c.ecp.fr}

\begin{abstract}
The scope of the present study is Eulerian modeling and simulation of polydisperse liquid sprays
undergoing droplet coalescence and evaporation.
The fundamental mathematical description is the
Williams spray equation governing the joint number density function $f (v , \bfu;\bfx,t)$
of droplet volume and velocity. Eulerian multi-fluid models have
already been rigorously derived from this equation in Laurent et al.~\cite{lmv04}.
The first key feature of the paper is the application of direct quadrature method of moments (DQMOM)
introduced by Marchisio and Fox \cite{marchisio05} to the Williams spray equation. Both the multi-fluid method and DQMOM yield systems of
Eulerian conservation equations with complicated interaction terms representing coalescence. In order to focus on the difficulties associated with treating size-dependent coalescence and to avoid numerical uncertainty issues associated with two-way coupling, only one-way coupling between the droplets and a given gas velocity field is considered. In order to validate and compare these approaches, the
  chosen configuration
is a self-similar 2D axisymmetrical decelerating nozzle with
sprays having various size distributions, ranging from smooth ones up to Dirac delta functions.
The second key feature of the paper is a thorough comparison
of the two approaches for various test-cases to a reference solution
obtained through a classical stochastic Lagrangian solver.
Both Eulerian models prove to
describe adequately spray coalescence and yield a very interesting alternative to the
Lagrangian solver.
The third key point of the study is a detailed description of the
limitations associated with each method, thus giving criteria
for their use as well as for their respective efficiency.
\end{abstract}
\end{frontmatter}

\section{Introduction}

In many industrial combustion applications such as Diesel engines,
fuel is stocked in condensed form and burned as a dispersed liquid phase
carried by a gaseous flow.
Two phase effects as well as the polydisperse character of the droplet size distribution (since the droplets dynamics depend on their inertia and are conditioned by size)
can significantly influence flame structure.
Size distribution effects are also encountered in a crucial way in
solid propellant rocket boosters, where the cloud of alumina particles
experiences coalescence and become polydisperse in size, thus determining their global dynamical behavior \cite{hylkema99,hylkema98}.
The coupling of dynamics, conditioned on particle size, with coalescence or aggregation as well as with evaporation can also be found in the study of fluidized beds \cite{tsuji98} and planet formation in solar nebulae \cite{bracco99,chavanis00}.
Consequently, it is important to have reliable models and numerical methods in order to be able to describe precisely the physics of two-phase flows where the dispersed phase is constituted of a cloud of particles of various sizes that can evaporate, coalesce or aggregate and also have their own inertia and size-conditioned dynamics.
Since our main area of interest is combustion, we will work with sprays throughout this paper, keeping in mind the broad application fields related to this study.

Generally speaking, two approaches for treating liquid sprays corresponding to two levels of description can be distinguished.
The first, associated with a full direct numerical simulation of the process, provides a model for the dynamics of the interface between the gas and liquid, as well as the exchanges of heat and mass between the two phases using various techniques such as the Volume Of Fluids (VOF) or Level Set methods \cite{aulisa03,herrmann05,josserand05,tanguy05}.
This ``microscopic" point of view is very rich in information on the detailed properties at a more local level concerning, for example, the resulting drag exerted on one droplet depending on its surroundings.
The second approach, based on a more global point of view, describes the droplets as a cloud of point particles for which the exchanges of mass, momentum and heat are described globally, using eventually correlations, and the details of the interface
behavior, angular momentum of droplets, detailed internal temperature distribution inside the droplet, etc., are not predicted.
Instead, a finite set of global properties such as mass, momentum, temperature are modeled.
Because it is the only one for which numerical simulations at the scale of a combustion chamber or in a free jet can be conducted, this ``mesoscopic" point of view will be adopted in the present paper.

Furthermore, we are interested in sprays where droplet interactions (e.g., coalescence) have to be taken into account, which corresponds to liquid volume fractions between $0.1\%$ and $1\%$. O'Rourke \cite{orourke81} classified the various regimes from the ``very thin spray'', which are transported by the gaseous carrier phase without influencing the gaseous phase, through the ``thin spray'' regime, for which there is two-way coupling through the momentum equation between the two phases, up to the 
``thick spray'' regime for which the volume fraction of liquid is high enough so that droplet-droplet interactions have to be taken into account, but is still low enough so that the liquid volume fraction is negligible as compared to the gaseous one. Because our primary focus is on the ability of Eulerian methods to capture droplet coalescence, our study is limited here to the ``thick spray'' regime.   By restricting our attention to one-way coupling, we can avoid difficulties (e.g., grid convergence) associated with using Lagrangian methods with two-way coupling, and it will thus be possible to make detailed comparisons between Eulerian and Lagrangian simulation results.

In the mesoscopic framework, there exists considerable interest in the development of numerical methods for simulating sprays \cite{hylkema98,hylkema99,miller99,miller00,lmv04,reveillon04}.
The principal physical processes that must be accounted for are (1) transport in real space, (2) droplet evaporation, (3) acceleration of droplets due to drag, and (4) coalescence of droplets leading to polydispersity.
The major challenge in numerical simulations is to account for the strong coupling between these processes.
Williams \cite{williams58} proposed a relatively simple transport equation based on kinetic theory that has proven to be a useful starting point for testing novel numerical methods for treating coalescing liquid sprays.
In the context of one-way coupling, the Lagrangian Monte-Carlo approach \cite{dukowicz80}, called Direct Simulation Monte-Carlo method (DSMC) by Bird \cite{bird94}, is generally considered to be more accurate than Eulerian methods for solving Williams equation.
However, its computational cost is high, especially in unsteady configurations. Moreover, in applications with two-way coupling, Lagrangian methods are difficult to couple accurately with Eulerian descriptions of the gas phase.  There is thus considerable impetus to develop Eulerian methods for describing sprays.
In this paper, we limit our attention to one-way coupling with a given (laminar) gas velocity field (i.e., one-way coupling with a given gas velocity field.)  Thus no turbulence models are required to close the spray equation.

In a recent paper Laurent et al.~\cite{lmv04} have demonstrated the capability of an Eulerian multi-fluid model to capture the physics of polydisperse evaporating sprays with one-way coupling.
This approach relies on the derivation of a semi-kinetic model from the Williams equation using a moment method for velocity, but keeping the continuous size distribution function.
This distribution function is then discretized using a ``finite-volume'' approach that yields conservation equations for mass, momentum (and eventually other properties such as temperature) of droplets in fixed size intervals called ``sections'' extending the original work of Tambour, Greenberg and collaborators \cite{greenberg86,greenberg93}.
Even though this approach has recently been extended to higher order by Laurent \cite{laurent06} and Dufour \cite{dufour_th,dufour05}, the necessity to discretize the size phase space can be a stumbling block in some applications.
Moment methods, on the other hand, do not encounter this limitation.

In this work, we apply the recently developed direct quadrature method of moments (DQMOM) \cite{marchisio05} to treat Williams equation in a Eulerian framework.
As its name implies, DQMOM is a moment method that closes the non-linear terms (e.g., droplet coalescence) using weighted quadrature points (abscissas) in phase space.
Such a closure relates to the construction of an approximated
number density function from a set of  moments under the form of a sum of Dirac delta functions,
the support of which corresponds to the abscissas. However, it is important 
to make a clear difference between such an Eulerian approach and the corresponding Lagrangian approach, for which the number density is approximated by a large number of numerical ``parcels''. The evolution of abscissas and the corresponding weights are governed by the dynamics of a few moments, whereas the evolution of the parcels are governed by the Williams equation since they are a stochastic discretization of this equation. Consequently, the DQMOM usually involves a very restricted number of unknowns on a Eulerian mesh, whereas the Lagrangian method involves a very large number of unknowns that are followed along their trajectories in phase space.

The DQMOM method distinguishes itself from other quadrature methods (e.g., QMOM \cite{mcgraw97,marchisio03}) by solving transport equations for the weights and abscissas directly (instead of transport equations for the moments).
The source terms for the transport equations depend on the physical processes involved.
For Williams equation, we show in Section~\ref{sec2} that laminar transport and drag result in source terms that are independent of the choice of moments and, in fact, are equivalent to those used in Lagrangian formulations.
When evaporation does not lead to the disappearance of droplets in finite time, this is also true for the evaporation process.
On the other hand, coalescence leads to a linear system for the source terms for which the coefficient matrix depends on the choice of moments.
The applicability of DQMOM to Williams equation thus depends on whether or not a particular choice(s) of moments can be found that leads to a non-singular linear system.  
When the evaporation law allows the disappearance of droplets in finite time the equations for the moments of the number density function not only involve unclosed integral terms, but also the flux of disappearing droplets, i.e.\ the pointwise value of the number density function at zero size. This quantity has then to be closed since it has a strong influence on the dynamics of the whole set of moments; it leads to a significant difficulty since it corresponds to the reconstruction of a pointwise value of the number density function from a set of its moments. 
In this study,  we propose a solution to this difficult issue.  
Note that because spatial transport is treated explicitly, it suffices  
to tackle the flux problem in the homogeneous case. 
We will see that a key point is to provide a flux closure that yields stable moment dynamics and a non-singular linear system in the DQMOM framework.

Let us also underline that the transport terms in the systems of conservation equations
for both Eulerian models are the same and given by pressureless gas dynamics. The structure of these transport terms
and the associated difficulties have been the subject of several studies and there are numerical methods designed
in order to treat the resulting singularities as shown in \cite{lmv04}. 
The question of the computational efficiency of
such Eulerian approaches (especially in coalescing systems) is a key question since these methods are intended to be used in
more realistic unsteady configurations as an alternative to the too costly Lagrangian methods for polydisperse sprays.
We have already studied this question in \cite{lmv04} where
the Eulerian multi-fluid approach was shown to offer a good precision with a relative low cost \cite{lmv04}.
Because of the similarity of the transport terms for both Eulerian approaches, the conclusions about the
computational efficiency presented in
\cite{lmv04} are also valid for the DQMOM method.  Consequently we
focus our study and comparisons on stationary configurations for which we are sure to have
a reference solution at our disposal and from which we can obtain firm conclusions about the capabilities of the various approaches.

In Section~\ref{test}, we present the chosen test configuration, which is a self-similar  2D axisymmetrical decelerating nozzle and sprays with two inlet distributions: a smooth monomodal function and Dirac delta functions.
We also discuss in detail the reasons (e.g., significant coalescence rates) for the choice of the test cases, and why they are particularly challenging for the various numerical methods.
Finally the Lagrangian solver, the numerical subtleties for obtaining the associated reference solution, as well as the multi-fluid method are then presented.
In Section~\ref{results}, we consider the results for the various test cases including  combinations of coalescence, linear evaporation in terms of volume (since it conserves the number of droplets and thereby eliminates the need to model the evaporative flux)
and the usual non-linear evaporation law (for which the evaporative flux must be modeled.)
We present results for the most difficult test cases,
designed to highlight the
challenges one would encounter in more realistic cases.
The results are compared to a reference solution obtained through a Lagrangian stochastic algorithm \cite{hylkema99}.
The advantages and limitations of the Eulerian methods are then analyzed in detail in terms of precision and efficiency.
It is shown that the DQMOM method offers very interesting features in a number of situations (e.g., strongly coalescing droplets), and is a good candidate for more complex configurations.

\section{DQMOM for Williams equation}\label{sec2}

The Williams transport equation \cite{williams58} for the joint volume, velocity number density function $f (v, \bfu; \bfx, t)$ is
\begin{equation}\label{weq}
\partial_{t} f + \bfu \cdot \partial_{\bfx} f
+ \partial_{v} \left( R_v f \right) + \partial_{\bfu} \cdot \left( \bfF f \right) = \Gamma,
\end{equation}
where $R_v$ is the evaporation rate, $\bfF$ is the drag force acting on the droplet, and $\Gamma$ is the coalescence term. Note that specific forms for the evaporation rate and drag law are not required for DQMOM. However, in this work we will consider one-way coupling with a given gas velocity that appears in $\bfF$.  Using standard assumptions \cite{lmv04}, we can write the coalescence term in two parts: $\Gamma = Q_{coll}^{-} + Q_{coll}^{+}$ where
\begin{align}
Q_{coll}^{-} &= - \int \int_{0}^{\infty} B(|\bfu - \bfu^*|, v, v^*)
f (v, \bfu) f (v^*, \bfu^*)  \sd v^* \sd \bfu^* , \label{qminus} \\
Q_{coll}^{+} &= \frac{1}{2} \int \int_{0}^{v} B(|\bfu^{\diamond} - \bfu^*|, v^{\diamond}, v^*)
f (v^{\diamond}, \bfu^{\diamond}) f (v^*, \bfu^*)  J \sd v^* \sd \bfu^*, \label{qplus}
\end{align}
$v^{\diamond} = v - v^*$, $\bfu^{\diamond} = (v \bfu - v^* \bfu^*)/(v - v^*)$, and $J = (v / v^{\diamond})^3$ is the Jacobian of the transform $(v, \bfu) \rightarrow (v^{\diamond}, \bfu^{\diamond})$ with fixed $(v^*, \bfu^*)$.  The collision frequency function $B$ is defined by
\begin{equation}\label{collf}
B(|\bfu - \bfu^*|, v, v^*) = E_{\rm coal}(|\bfu - \bfu^*|, v, v^*) \beta(v, v^*) |\bfu - \bfu^*|,
\end{equation}
where $E_{\rm coal}$ is
the coalescence efficiency probability, which, based upon the size of droplets and
the relative velocity, discriminates between rebound and coalescence,
and
\begin{equation}
\beta (v , v^*) = \pi \left[ \left( \frac{3 v}{4 \pi} \right)^{1/3} + \left( \frac{3 v^*}{4 \pi} \right)^{1/3} \right]^2 .
\end{equation}
For simplicity, we will take $E_{\rm coal} =0$ (no coalescence) or $E_{\rm coal} =1$; however, any other functional form could be used in the derivation that follows.
A more general version of the spray equation would include the droplet temperature and molecular composition.
For simplicity, we consider only the volume and velocity in this work.
Finally note that adding spatial diffusion terms in Eq.~(\ref{weq}) would generate additional terms in DQMOM \cite{marchisio05}.

One of the principal mathematical difficulties when developing Eulerian solvers for Eq.~(\ref{weq}) is the accurate treatment of the coalescence term.  Indeed, the integral form of $\Gamma$ leads to highly non-local and non-linear interactions in volume-velocity phase space.  A ``direct'' Eulerian solver would require discretization of the high-dimensional phase space (in addition to real space), and would thus be computationally intractable.  In contrast, multi-fluid models discretize only the volume phase space and use the average velocity conditioned on droplet size (i.e., the mono-kinetic assumption \cite{laurent01}), while moment methods (such as DQMOM) provide closures based on a finite set of moments. 
Before applying DQMOM to Eq.~(\ref{weq}), we should note that the coalescence term is defined such that the moments representing mass and momentum are conserved:
\begin{equation}\label{cmass}
\int  \rho v \Gamma (v, \bfu) \sd v \sd \bfu = 0
\end{equation}
and
\begin{equation}\label{cmom}
\int \rho v \bfu \Gamma (v, \bfu) \sd v \sd \bfu = 0,
\end{equation}
where the liquid density $\rho$ is assumed to be constant.
These conservation properties must be retained in numerical approximations used to treat Eq.~(\ref{weq}) (as we shall see is the case with DQMOM).

The DQMOM approximates the density function by weighted delta functions in volume-velocity phase space \cite{fox03,marchisio05}:
\begin{equation}\label{dqmom1}
f (v, \bfu ) = \sum_{n=1}^{N} w_n \delta (v - v_n) \delta ( \bfu - \bfu_n )
\end{equation}
where $\delta ( \bfu - \bfu_n ) \equiv \delta ( u_1 - u_{1,n} ) \delta ( u_2 - u_{2,n} ) \delta ( u_3 - u_{3,n} )$.
Note that in this formulation, the weights $w_n$ and abscissas ($v_n$, $\bfu_n$) are Eulerian fields.
Application of DQMOM results in closed transport equations for the number density, mass density, and momentum density, respectively, of the form:
\begin{gather}
 \partial_{t} w_n + \partial_{\bfx} \cdot \left( w_n \bfu_n \right)
= a_n , \label{number}   \\
 \partial_{t} \left( w_n \rho v_n \right) +  \partial_{\bfx} \cdot \left( w_n \rho v_n \bfu_n \right)
= \rho b_n , \label{mass}  \\
\intertext{and}
 \partial_{t} \left( w_n \rho v_n \bfu_n \right) + \partial_{\bfx} \cdot \left( w_n \rho v_n \bfu_n \bfu_n \right)
= \rho \bfc_n , \label{momentum}
\end{gather}
where $a_n$, $b_n$, and $\bfc_n$ are source terms that are found from the right-hand side of Eq.~(\ref{weq}) as described below.
These equations can be solved with appropriate initial and boundary conditions to find the fields $w_n (\bfx, t)$ and ($v_n (\bfx, t)$, $\bfu_n (\bfx, t)$) appearing in Eq.~(\ref{dqmom1}).
Note that Eqs.~(\ref{number}--\ref{momentum}) are equivalent to an Eulerian multi-fluid model \cite{lmv04}, but with the source terms on the right-hand side determined using DQMOM.

The DQMOM approximation for the moments of the number density function are found directly from Eq.~(\ref{dqmom1}):
\begin{equation}\label{dqmom2}
\langle v^k u_1^l u_2^m u_3^p \rangle \equiv
\int v^{k} u_{1}^{l} u_2^{m} u_3^{p}
f (v, \bfu ) \sd v \sd \bfu
= \sum_{n=1}^{N} w_n v_n^{k} u_{1, n}^{l} u_{2, n}^{m} u_{3, n}^{p} .
\end{equation}
The fundamental idea behind DQMOM is that we should choose the weights and abscissas such that as many moments as possible are determined by the moment transport equations found from Eq.~(\ref{weq}).
Note that there are a total of $N$ weights, $N$ volume abscissas, and $3N$ velocity abscissas and (equivalently) $5N$ unknown source terms in Eqs.~(\ref{number}--\ref{momentum}).
We will thus need to choose $5N$ independent moments to determine the source terms. We will return to the subject of how to choose the moments in Section~\ref{decomp}. The procedure for using these moments to find the source terms is described next.

\subsection{Space and time derivatives}

The space and time derivatives in Eq.~(\ref{weq}) generate the corresponding terms in Eqs.~(\ref{number}--\ref{momentum}).
These are found by formally inserting Eq.~(\ref{dqmom1}), and differentiating:
\begin{multline}\label{stder}
\partial_{t} f + \bfu \cdot \partial_{\bfx} f =
\sum_{n=1}^{N}
\delta (v - v_n) \delta ( \bfu - \bfu_n ) \left[ \partial_{t} w_n +  \partial_{\bfx} \cdot ( \bfu_n w_n ) \right] \\
- \sum_{n=1}^{N} w_n \delta^{(1)} (v - v_n) \delta ( \bfu - \bfu_n ) \left[ \partial_{t} v_n +  \bfu_n \cdot \partial_{\bfx}  v_n  \right] \\
- \sum_{n=1}^{N} w_n \delta (v - v_n) \bfdelta^{(1)} ( \bfu - \bfu_n ) \cdot \left[ \partial_{t} \bfu_n +  \bfu_n \cdot \partial_{\bfx}  \bfu_n  \right]
\end{multline}
where $\delta^{(1)} ( \psi ) = \rd \delta ( \psi ) / \rd \psi$ and $\bfdelta^{(1)} ( \bfpsi )$ is a vector with components $\delta^{(1)}_{1} ( \bfpsi ) = \delta^{(1)} ( \psi_1 ) \delta ( \psi_2 ) \delta ( \psi_3 )$, $\delta^{(1)}_{2} ( \bfpsi ) = \delta ( \psi_1 ) \delta^{(1)} ( \psi_2 ) \delta ( \psi_3 )$, and $\delta^{(1)}_{3} ( \bfpsi ) = \delta ( \psi_1 ) \delta ( \psi_2 ) \delta^{(1)} ( \psi_3 )$.
Using the definitions of the source terms, Eq.~(\ref{stder}) can be rewritten as
\begin{multline}\label{stder1}
\partial_{t} f + \bfu \cdot \partial_{\bfx} f =
\sum_{n=1}^{N} \left[ \delta (v - v_n) \delta ( \bfu - \bfu_n ) + v_n \delta^{(1)} (v - v_n) \delta ( \bfu - \bfu_n ) \right] a_n \\
- \sum_{n=1}^{N} \left[ \delta^{(1)} (v - v_n) \delta ( \bfu - \bfu_n ) - v_n^{-1} \delta (v - v_n) \bfdelta^{(1)} ( \bfu - \bfu_n ) \cdot \bfu_n \right] b_n \\
- \sum_{n=1}^{N} v_n^{-1} \delta (v - v_n) \bfdelta^{(1)} ( \bfu - \bfu_n ) \cdot \bfc_n .
\end{multline}
Note that this expression is linear in the source terms ($a_n$, $b_n$, $\bfc_n$).

The next step is to apply the moment transform to Eq.~(\ref{stder1}).
Formally, this yields
\begin{multline}\label{momt1}
\int v^{k} u_{1}^{l} u_2^{m} u_3^{p}
\left( \partial_{t} f + \bfu \cdot \partial_{\bfx} f \right) \sd v \sd \bfu
=
\sum_{n=1}^{N} ( 1 - k ) v_n^{k} u_{1,n}^{l} u_{2,n}^{m} u_{3,n}^{p} a_n \\
+ \sum_{n=1}^{N} (k - l - m - p) v_n^{k-1} u_{1,n}^{l} u_{2,n}^{m} u_{3,n}^{p} b_n \\
+ \sum_{n=1}^{N} v_n^{k-1} u_{1,n}^{l} u_{2,n}^{m} u_{3,n}^{p}
\left( l u_{1,n}^{-1} c_{1,n} + m u_{2,n}^{-1} c_{2,n} + p u_{3,n}^{-1} c_{3,n} \right),
\end{multline}
where, unless otherwise noted, the definite integrals cover all of phase space.
The next step is to consider the terms in Eq.~(\ref{weq}) that correspond to transport in volume-velocity phase space.

\subsection{Phase-space transport}

We begin by rewriting Eq.~(\ref{weq}) as
\begin{equation}\label{weq1}
\partial_{t} f + \bfu \cdot \partial_{\bfx} f  = P ,
\end{equation}
where the phase-space transport terms are defined by
\begin{equation}\label{phase}
P \equiv - \partial_{v} \left( R_v f \right)
- \partial_{\bfu} \cdot \left( \bfF f \right) + \Gamma.
\end{equation}
We can then define the moment transform of the phase-space terms by
\begin{equation}\label{dqmom3}
P (k,l,m,p) \equiv
\int v^{k} u_{1}^{l} u_2^{m} u_3^{p} P \sd v \sd \bfu .
\end{equation}
Note that if the moments $P (k,l,m,p)$ are known, Eq.~(\ref{momt1}) forms a linear system that can be solved to find the unknown source terms.
We can compute the phase-space moments using Eq.~(\ref{phase}):
\begin{equation}
P (k,l,m,p) =
- \int v^{k} u_{1}^{l} u_2^{m} u_3^{p}
\left[ \partial_{v} \left( R_v f \right)
+ \partial_{\bfu} \cdot \left( \bfF f \right) - \Gamma \right]
\sd v \sd \bfu . \label{dqmom4}
\end{equation}
As shown next, the integrals on the right-hand side can be expressed in terms of the weights and abscissas, and a flux term corresponding to disappearance of droplets due to evaporation.

Starting with the evaporation term in Eq.~(\ref{dqmom4}), we can use integration by parts to find
\begin{equation}\label{evp1}
\int_{0}^{\infty} v^{k} \partial_{v} \left( R_v f \right) \sd v
=
-
\delta_{k0} R_v (0 , \bfu ) f (0 , \bfu)
-
\int_{0}^{\infty} k v^{k-1} R_v (v , \bfu ) f \sd v ,
\end{equation}
where $\delta_{k0}$ is the Kronecker delta.  Using Eq.~(\ref{dqmom1}) in the final integral, we find
\begin{multline}\label{evp2}
\int v^{k} u_{1}^{l} u_2^{m} u_3^{p} \partial_{v} \left( R_v f \right) \sd v \sd \bfu
=
-
\delta_{k0} \psi(t) u_{f1}^{l} u_{f2}^{m} u_{f3}^{p} \\
-
\sum_{n=1}^{N} k w_n v_n^{k-1} u_{1,n}^{l} u_{2,n}^{m} u_{3,n}^{p} R_v (v_n , \bfu_n ),
\end{multline}
where $\psi (t)$ is the evaporative flux of droplets at zero size 
and $\bfu_{f}$ is the velocity of droplets with zero volume (which will normally correspond to the fluid velocity).
Note that the first term on the right-hand side of Eq.~(\ref{evp2}) will be non-zero only for $k=0$, and corresponds to the loss of droplets due to evaporation.
A fundamental question when applying DQMOM to evaporation problems is how to determine $\psi(t)$ from the weights and abscissas.
The value of $\psi(t)$ corresponds to the value of the number density function
at zero size, and in the case of the $d^2$ evaporation law, it is precisely the 
value of the number density as a function of droplet surface, which has no reason to be zero in general. Determining the value of $\psi(t)$, a pointwise information, from the values of moments is clearly a difficult task, for which we will propose a solution in the next subsection. 
On the other hand, the second term on the right-hand side of Eq.~(\ref{evp2}) is non-zero only for $k>0$, and appears in closed form.

Turning next to the drag-force term in Eq.~(\ref{dqmom4}), we can use integration by parts to find
\begin{equation}
\int u_j^l \partial_{u_j} \left( F_j f \right) \sd \bfu
=
- \int l u_j^{l-1} F_j f  \sd \bfu
\quad \text{for $j = 1,2,3$}. \label{drag2}
\end{equation}
Thus, the drag-force term becomes
\begin{multline}
\int v^{k} u_{1}^{l} u_2^{m} u_3^{p}
\partial_{\bfu} \cdot \left( \bfF f \right) \sd v \sd \bfu = \\
-
\sum_{n=1}^{N} w_n v_n^k u_{1,n}^{l} u_{2,n}^{m} u_{3,n}^{p} \\
\left[ l u_{1,n}^{-1} F_1 (v_n , \bfu_n ) +  m u_{2,n}^{-1} F_2 (v_n , \bfu_n )
    +  p u_{3,n}^{-1} F_3 (v_n , \bfu_n ) \right]. \label{du1}
\end{multline}
Note that this term appears in closed form.

Turning now to the coalescence term, we will treat each of the two parts $Q^{-}_{coll}$ and $Q^{+}_{coll}$ separately.
The first part yields in a straightforward manner
\begin{equation}
\int v^{k} u_{1}^{l} u_2^{m} u_3^{p}
Q^{-}_{coll} \sd v \sd \bfu =
-
\sum_{n=1}^{N} \sum_{q=1}^{N} w_n w_q v_{n}^k u_{1,n}^{l} u_{2,n}^{m} u_{3,n}^{p} B ( |\bfu_n - \bfu_q |, v_n , v_q )
. \label{cmu3}
\end{equation}
The second part requires a change in the order of integration, and a change of variables:
\begin{multline}\label{collcv}
\int \left( \int_{0}^{v}
h(v, \bfu) B(|\bfu^{\diamond} - \bfu^*|, v^{\diamond}, v^*)
f (v^{\diamond}, \bfu^{\diamond}) f (v^*, \bfu^*)
J \sd v^* \right) \sd v \sd \bfu^* \sd \bfu \\
=
\int \left( \int_{v^*}^{\infty}
h(v, \bfu) B(|\bfu^{\diamond} - \bfu^*|, v^{\diamond}, v^*)
f (v^{\diamond}, \bfu^{\diamond}) f (v^*, \bfu^*)
J \sd v \right) \sd v^* \sd \bfu^* \sd \bfu \\
=
\int
h \left( v^* + v^\diamond , \frac{v^* \bfu^* + v^\diamond \bfu^\diamond }{v^* + v^\diamond} \right) \\
B(|\bfu^{\diamond} - \bfu^*|, v^{\diamond}, v^*)
f (v^{\diamond}, \bfu^{\diamond}) f (v^*, \bfu^*)  \sd v^*
\sd v^\diamond \sd \bfu^* \sd \bfu^\diamond ,
\end{multline}
where $h$ is an arbitrary function of $v$ and $\bfu$. It then follows that
\begin{multline}
\int v^{k} u_{1}^{l} u_2^{m} u_3^{p}
Q^{+}_{coll} \sd v \sd \bfu =
\frac{1}{2}
\sum_{n=1}^{N} \sum_{q=1}^{N} w_n w_q ( v_n + v_q )^k
\left( \frac{v_n u_{1,n} + v_q u_{1,q}}{v_n + v_q} \right)^{l} \\ \times
\left( \frac{v_n u_{2,n} + v_q u_{2,q}}{v_n + v_q} \right)^{m}
\left( \frac{v_n u_{3,n} + v_q u_{3,q}}{v_n + v_q} \right)^{p}
B ( |\bfu_n - \bfu_q |, v_n , v_q )
. \label{cpu3}
\end{multline}
Note that the right-hand side of this expression is in closed form.

Collecting together all of the terms, the moments appearing on the right-hand sides of Eqs.~(\ref{evp2}--\ref{cpu3}) become
\begin{multline}
P (k,l,m,p) =
\delta_{k0} \psi(t) u_{f1}^{l} u_{f2}^{m} u_{f3}^{p}
+
\sum_{n=1}^{N} k w_n v_n^{k-1} u_{1,n}^{l} u_{2,n}^{m} u_{3,n}^{p} R_v (v_n , \bfu_n ) \\
+
\sum_{n=1}^{N} w_n v_n^k u_{1,n}^{l} u_{2,n}^{m} u_{3,n}^{p}
\left[ l u_{1,n}^{-1} F_1 (v_n , \bfu_n )
+  m u_{2,n}^{-1} F_2 (v_n , \bfu_n )
+  p u_{3,n}^{-1} F_3 (v_n , \bfu_n ) \right] \\
+
\frac{1}{2}
\sum_{n=1}^{N} \sum_{q=1}^{N} w_n w_q \\
\left[ ( v_n + v_q )^k
\left( \frac{v_n u_{1,n} + v_q u_{1,q}}{v_n + v_q} \right)^{l}
\left( \frac{v_n u_{2,n} + v_q u_{2,q}}{v_n + v_q} \right)^{m}
\left( \frac{v_n u_{3,n} + v_q u_{3,q}}{v_n + v_q} \right)^{p}  \right. \\ \left.
\vphantom{\left( \frac{v_n u_{1,n} + v_q u_{1,q}}{v_n + v_q} \right)^{l}}
- v_{n}^k u_{1,n}^{l} u_{2,n}^{m} u_{3,n}^{p}
- v_{q}^k u_{1,q}^{l} u_{2,q}^{m} u_{3,q}^{p}
\right]
B ( |\bfu_n - \bfu_q |, v_n , v_q ). \label{volf}
\end{multline}
Note that due to the form of the coalescence term, the moments conserve mass ($P (1,0,0,0) =0$) and momentum ($P (1,1,0,0) = P (1,0,1,0) = P (1,0,0,1) = 0$) when evaporation and drag are null.  Thus, the weights and abscissas in the DQMOM representation will keep the same conservation properties as the original model (i.e., as Eq.~(\ref{weq})).

Comparing the terms in Eqs.~(\ref{momt1}) and (\ref{volf}), we can note that the evaporation and drag terms in the DQMOM representation can be solved for explicitly.  Thus, the source terms can be written as
\begin{align}
b_n &= b_n^* + w_n R_v (v_n , \bfu_n ), \\
\bfc_{n} &= \bfc_{n}^* + w_n \bfu_{n} R_v (v_n , \bfu_n ) + w_n v_n \bfF (v_n , \bfu_n ),
 \label{bcs}
\end{align}
where source terms $a_n$, $b_n^*$ and $\bfc_{n}^*$ in the transport equations are found by solving the linear system
\begin{multline}
 \sum_{n=1}^{N} ( 1 - k ) v_n^{k} u_{1,n}^{l} u_{2,n}^{m} u_{3,n}^{p} a_n
+ \sum_{n=1}^{N} (k - l - m - p) v_n^{k-1} u_{1,n}^{l} u_{2,n}^{m} u_{3,n}^{p} b_n^* \\
+ \sum_{n=1}^{N} v_n^{k-1} u_{1,n}^{l} u_{2,n}^{m} u_{3,n}^{p}
\left( l u_{1,n}^{-1} c_{1,n}^* + m u_{2,n}^{-1} c_{2,n}^* + p u_{3,n}^{-1} c_{3,n}^* \right)
= P^* (k,l,m,p), \label{u3ns}
\end{multline}
with the right-hand side given by
\begin{multline}
P^* (k,l,m,p) = - \delta_{k0} \psi u_{f1}^{l} u_{f2}^{m} u_{f3}^{p}
+
\frac{1}{2}
\sum_{n=1}^{N} \sum_{q=1}^{N} w_n w_q \\
\left[ ( v_n + v_q )^k
\left( \frac{v_n u_{1,n} + v_q u_{1,q}}{v_n + v_q} \right)^{l}
\left( \frac{v_n u_{2,n} + v_q u_{2,q}}{v_n + v_q} \right)^{m}
\left( \frac{v_n u_{3,n} + v_q u_{3,q}}{v_n + v_q} \right)^{p} \right. \\ \left.
\vphantom{\left( \frac{v_n u_{1,n} + v_q u_{1,q}}{v_n + v_q} \right)^{l}}
- v_{n}^k u_{1,n}^{l} u_{2,n}^{m} u_{3,n}^{p}
- v_{q}^k u_{1,q}^{l} u_{2,q}^{m} u_{3,q}^{p}
\right]
B ( |\bfu_n - \bfu_q |, v_n , v_q ) . \label{volfs}
\end{multline}
The expression for the source terms (Eq.~\ref{u3ns}) completes the derivation of the DQMOM transport equations for the Williams spray equation.

In the absence of coalescence, Eq.~(\ref{volfs}) is particularly simple.
Thus, the pure evaporation case is an interesting limit case for which $a_n$, $b_n^*$, and $\bfc_{n}^*$ will be non-zero only if the evaporative flux $\psi$ is non-zero.
However, the evaporative flux cannot be determined by moment constraints alone (see Section~\ref{ef}).
If the evaporative flux is assumed to be null, the zero-order moment will remain unchanged in the absence of coalescence
as long as some abscissa crosses the zero size limit and yields a pointwise 
singular and infinite flux as in Lagrangian methods when some parcels reach the zero size limit. However, as mentioned in the Introduction,  since there are only a few abscissas that describe the moment dynamics, 
such a singular behavior is not ideal for smooth number density functions (whereas it is the correct one if the number density function is a sum of Dirac delta function from the beginning as in the bimodal case that will be studied later). Consequently we need an evaluation of this flux function that guarantees
a smooth flux as a function of time for smooth distribution functions.
Even when coalescence is included, the moments may be poorly estimated if the evaporative flux is neglected.
An example of such behavior can be found in the work of Mossa \cite{mossa05_th} where the droplet size distribution is presumed to be log-normal and where the evaporative flux at zero size is neglected, leading to numerical difficulties and a poor prediction of the second moment.
 Thus, we will use a separate procedure, described next, to approximate the contribution due to the evaporative flux 
that yields a continuous in time flux, as well as a guarantee that the abscissas never cross the zero size limit.

\subsection{Evaporative flux}\label{ef}

The source terms cannot be computed directly from the moment constraints in Eq.~(\ref{volfs}) because the evaporative flux is unknown.
We must therefore apply additional (or different) constraints to determine all of the unknowns. Considering only evaporation
and setting drag and coalescence to zero in the right-hand side of Eq.~(\ref{u3ns}), we obtain the following linear system:
\begin{multline}
 \sum_{n=1}^{N} ( 1 - k ) v_n^{k} u_{1,n}^{l} u_{2,n}^{m} u_{3,n}^{p} a_n
+ \sum_{n=1}^{N} (k - l - m - p) v_n^{k-1} u_{1,n}^{l} u_{2,n}^{m} u_{3,n}^{p} b_n^* \\
+ \sum_{n=1}^{N} v_n^{k-1} u_{1,n}^{l} u_{2,n}^{m} u_{3,n}^{p}
\left( l u_{1,n}^{-1} c_{1,n}^* + m u_{2,n}^{-1} c_{2,n}^* + p u_{3,n}^{-1} c_{3,n}^* \right) \\
+
\delta_{k0} u_{f1}^{l} u_{f2}^{m} u_{f3}^{p} \psi = 0 \label{evap1}
\end{multline}
with $5N+1$ unknowns $a_n$, $b_n^*$, $\bfc_{n}^*$ and $\psi$.
Note that because the right-hand side is null, only trivial solutions can be found using moment constraints.
We will therefore introduce ratio constraints of the form
\begin{equation*}
 \frac{\rD}{\rD t}
\left( \frac{w_n}{w_{n+1}} \right)_{\rm evap} =0 ,
\quad
\frac{\rD}{\rD t}
\left( \frac{v_n}{v_{n+1}} \right)_{\rm evap} = 0
\quad \text{and} \quad
\frac{\rD}{\rD t}
\left( \frac{u_{j  n}}{u_{j  n+1}} \right)_{\rm evap} =0,
\end{equation*}
which are applied only for the changes due to evaporation.
These constraints are motivated by the behavior of the weights and abscissas corresponding to sufficiently smooth and continuous density functions.
For example, if the surface density function is exponential and the evaporation rate is proportional to the surface area of a droplet, then the abscissas remain constant and the weights decrease monotonely.
On the other hand, for 
singular density functions (e.g., composed of delta functions), the ratio constraints are expected to perform poorly.
We will look more closely at this issue in Section~\ref{results}.

It can be observed that the choice of $k=l=m=p=0$ in Eq.~(\ref{evap1}) leads to
\begin{equation}\label{connum}
\psi = - \sum_{n=1}^{N} a_n .
\end{equation}
Thus, the evaporative flux depends only on $a_n$.
Note that physically $\psi \ge 0$.
Hence, if the value computed for $\psi$ from Eq.~(\ref{connum}) is negative (which is possible for very general evaporation rates), then $a_n$, $b_n^*$, $\bfc_{n}^*$ and $\psi$ are set equal to zero.
However, for the evaporation rate considered in this paper, it can be shown that the flux will be non-negative.

Conservation of mass ($k=1$ and $l=m=p=0$ in Eq.~(\ref{evap1})) leads to
\begin{equation}\label{conmass}
\sum_n b_n^* =0.
\end{equation}
Applying the ratio constraint for the abscissas yields
\begin{equation}\label{rca}
w_{n+1} v_{n+1} b^*_n - w_n v_n b^*_{n+1} = E_n
\quad \text{for $n=1,\ldots, N-1$};
\end{equation}
where the right-hand side is defined by
\begin{equation}\label{dn}
E_n =
w_n w_{n+1} \left[ v_{n} R_v (v_{n+1}) - v_{n+1} R_v (v_{n}) \right].
\end{equation}
Note that in order for there to be an evaporative flux, we will normally have $E_n \ge 0$ for all $n$ (assuming that $v_1<v_2<\ldots<v_N$).
The case where $E_n=0$ occurs when $R_v(v)$ is proportional to $-v$ (i.e., the evaporation rate is proportional to the droplet volume).
The more common case where $E_n > 0$ occurs when $R_v(v)$ is proportional to $-v^{1/3}$ (i.e., the droplet surface area decreases linearly).
In general, $R_v(v) \propto -v^{\gamma}$ with $\gamma < 1$ will lead to positive $E_n$.
The physical interpretation for this difference is that for $\gamma < 1$ the droplets will disappear due to evaporation in a finite time, while for $\gamma \ge 1$ the disappearance time is infinite.
The linear system formed from Eqs.~(\ref{conmass}) and (\ref{rca}) can be solved separately to find $b_n^*$.

Conservation of momentum ($k=1$ and $l$, $m$, or $p=1$ in Eq.~(\ref{evap1})) leads to
\begin{equation}\label{conmom}
 \sum_{n=1}^{N} \bfc_{n}^* = \bfzero .
\end{equation}
Likewise, the ratio constraint for each component of the velocity yields
\begin{equation}\label{rcv}
 w_{n+1} v_{n+1} u_{j  n+1} c^*_{j  n}
- w_n v_n u_{j  n} c^*_{j  n+1} =
u_{j  n} u_{j  n+1}  E_n
\quad \text{for $n=1,\ldots, N-1$}.
\end{equation}
Together with Eq.~(\ref{conmom}), this equation can be solved separately for each component ($j= 1,2,3$) to find $\bfc_{n}^*$.

The ratio constraint for the weights yields $N-1$ equations for $a_n$:
\begin{equation}\label{rcw}
 w_{n+1} a_n - w_n a_{n+1} = 0 \quad \text{for $n=1, \ldots, N-1$}.
\end{equation}
Note that this constraint is satisfied by $a_n = \alpha w_n$ where $\alpha$ is unknown. We must therefore choose one independent moment in Eq.~(\ref{evap1}) in order to solve for $\alpha$.  Since $b_n^*$ and $\bfc_{n}^*$ are known, we can rearrange Eq.~(\ref{evap1}) as
\begin{multline}
\alpha \sum_{n=1}^{N}
\left[ ( k - 1 ) v_n^{k} u_{1,n}^{l} u_{2,n}^{m} u_{3,n}^{p}
+
\delta_{k0} u_{f1}^{l} u_{f2}^{m} u_{f3}^{p} \right] w_n
= \\
 \sum_{n=1}^{N} (k - l - m - p) v_n^{k-1} u_{1,n}^{l} u_{2,n}^{m} u_{3,n}^{p} b_n^* \\
+ \sum_{n=1}^{N} v_n^{k-1} u_{1,n}^{l} u_{2,n}^{m} u_{3,n}^{p}
\left( l u_{1,n}^{-1} c_{1,n}^* + m u_{2,n}^{-1} c_{2,n}^* + p u_{3,n}^{-1} c_{3,n}^* \right), \label{evap2}
\end{multline}
which can be solved with $k \neq 1$ to find $\alpha$.
If we choose, for example, $k=2$ and $l=m=p=0$ as the independent moment, then the constraint becomes
\begin{equation}
 \alpha  =  2 \sum_{n=1}^{N} v_n b_n^*
\left/ \sum_{n=1}^{N} v_n^{2} w_n \right.
\label{evap3a}
\end{equation}
and $\alpha$ depends only on $b_n^*$.
However, if we choose $k=2$ and $l=m=p=1$, then the constraint becomes
\begin{multline}
 \alpha  =
 \sum_{n=1}^{N} v_n
\left(
  u_{2  n} u_{3 n} c_{1 n}^*
+ u_{1  n} u_{3 n} c_{2 n}^*
+ u_{1  n} u_{2 n} c_{3 n}^*
- u_{1  n} u_{2 n} u_{3 n} b_n^*
\right) \\
\left/ \sum_{n=1}^{N} v_n^{2} u_{1 n} u_{2 n} u_{3 n} w_n \right. .
\label{evap3}
\end{multline}
For this choice, $\alpha$ is independent of $\bfu_{f}$.
A choice that leads to a fully coupled system is $k=2$, $l=2$, $m=p=0$, which yields
\begin{equation}
 \alpha
=
2\sum_{n=1}^{N} v_n u_{1  n} c_{1 n}^*
\left/ \sum_{n=1}^{N} v_n^2 u_{1,n}^2 w_n \right.
\label{evap3b}
\end{equation}
or $k=m=p=0$ and $l=1$, which yields
\begin{equation}
 \alpha
=
\sum_{n=1}^{N} v_n^{-1}
\left( u_{1,n} b_n^* - c_{1,n}^* \right)
\left/ \sum_{n=1}^{N} \left( u_{1,n} - u_{f1} \right) w_n \right.
. \label{evap4}
\end{equation}
Note that when $v_n \rightarrow 0$, we have $u_{1,n} \rightarrow u_{f1}$ and $c_{1,n}^* \rightarrow u_{f1} b_n^*$; hence, this last constraint is consistent with this limiting behavior.
These choices are asymmetric in the velocity components, and thus do not treat all components the same.
A ``symmetric'' choice with similar properties is $k=2$ and $l=m=p=2$ or $k=0$ and $l=m=p=1$, which lead to a more complicated constraint.  The ``best'' choice will most likely be problem dependent.
In our test cases, the choices with $k=2$ give similar results, better than the ones with $k=0$.
The calculations are thus done with the value of $\alpha$ given in Eq.~(\ref{evap3a}):
this value is the simplest and can be shown to be non-positive as soon as $E_n\ge 0$, at least for the case $N=2$.

In summary, the contribution due to evaporation is estimated by first solving separate linear systems for $b_n^*$ and $\bfc_{n}^*$.
The estimate for $a_n = \alpha w_n$ is found using an independent moment constraint from Eq.~(\ref{evap2}) to find $\alpha$.
Finally, the evaporative flux $\psi$ is computed from Eq.~(\ref{connum}), and should be non-negative.
If $\psi$ is negative (or equivalently if $\alpha$ is positive), then the contribution due to evaporation is null.
The contribution due to coalescence is found by solving a linear system of the form of Eq.~(\ref{u3ns}) where the right-hand side is given by Eq.~(\ref{volfs}) with $\psi=0$.
 As described below, the final source terms ($a_n$, $b_n^*$, $\bfc^*_{n}$) are found simply by adding together the contributions from the evaporative flux and coalescence.

\subsection{DQMOM linear system}\label{decomp}

The DQMOM representation of Williams spray equation is given by the transport equations for the weights and abscissas (Eqs.~(\ref{number}--\ref{momentum})).
The source terms for these equations are found by solving the linear system as described above.
The exact form of the linear system depends on the choice of moments.
This choice, in turn, will determine if the system is well defined in the sense that the coefficient matrix is non-singular.  
A choice of moments that is consistent with the mono-kinetic assumption used in the multi-fluid model is to consider only moments of orders zero and one in the velocity components (i.e., $l,m,p \in \{0,1 \}$).  In this work, in order to make direct comparisons with the multi-fluid model, we will limit our consideration to such moments.  In general, this choice of moments should allow for the best possible description of droplet coalescence, while at the same time ensuring that droplet mass and momentum are conserved.

A choice of $5N$ moments that has been found to be non-singular is
\begin{align}
\nonumber
& k=(i-1)/3 \quad i\in \{1,\ldots,2N\} \quad &\mbox{with} \quad &l=m=p=0 \\
\nonumber
& k=i \quad i\in \{1,\ldots,N\} \quad &\mbox{with}  \quad &l=1, \quad m=p=0\\
\nonumber
& k=i \quad i\in \{1,\ldots,N\} \quad &\mbox{with} \quad &m=1, \quad l=p=0\\
& k=i \quad i\in \{1,\ldots,N\} \quad &\mbox{with} \quad &p=1, \quad l=m=0.
\label{set3i}
\end{align}
For $N\ge 2$, this choice of moments includes the surface area and the volume of the droplets, which are important variables for evaporating spray, as well as their momentum.
The linear system can then be written in matrix form (showing only non-zero components) as
\begin{equation}\label{linsys}
\begin{bmatrix}
 \bfA_1 & \bfA_2 & \bfE_1 & \bfE_2 & \bfE_3 \\
 \bfA_3 & \bfA_4 &        &        &        \\
 \bfB_1 & \bfC_1 & \bfD_1 &        &        \\
 \bfB_2 & \bfC_2 &        & \bfD_2 &        \\
 \bfB_3 & \bfC_3 &        &        & \bfD_3
\end{bmatrix}
\begin{bmatrix}
\bfa   \\
\bfb^*   \\
\bfc_1^* \\
\bfc_2^* \\
\bfc_3^*
\end{bmatrix}
=
\begin{bmatrix}
\bfP_a \\
\bfP_b \\
\bfP_1 \\
\bfP_2 \\
\bfP_3
\end{bmatrix},
\end{equation}
where the matrices $\bfA_j$, $\bfB_j$, $\bfC_j$, $\bfD_{j}$ and $\bfE_{j}$ are all $N$$\times$$N$, and $\bfa$, $\bfb^*$ and $\bfc_j^*$ are column vectors formed from the components $a_n$, $b_n^*$ and $c_{j n}^*$, respectively.
In general, the exact definitions of the other matrices will depend on which constraints are used to define the system, i.e., Eq.~(\ref{u3ns}) or those described in Section~\ref{ef}.
Nevertheless, the form of the linear system is the same in all cases.
As noted earlier, the linear system is solved twice at each time step.
First with the matrices for the evaporative flux without coalescence (i.e., $\bfA_3 = \bfB_j = \bfC_j = \bfzero$), and second with the matrices for coalescence without evaporation (i.e., $\bfE_j = \bfzero$).
The unknowns $\bfa$, $\ldots$, $\bfc_3^*$ are found by adding the two solutions.

As discussed earlier, for the evaporative step the linear system can be decomposed into five $N$$\times$$N$ systems that can be solved sequentially.
Likewise, for the coalescence step $\bfa$ and $\bfb^*$ can be found separately by solving a $2N$$\times$$2N$ system:
\begin{equation}\label{linsys2}
\begin{bmatrix}
 \bfA_1 & \bfA_2 \\
 \bfA_3 & \bfA_4
\end{bmatrix}
\begin{bmatrix}
\bfa   \\
\bfb^*
\end{bmatrix}
=
\begin{bmatrix}
\bfP_a \\
\bfP_b
\end{bmatrix},
\end{equation}
and then each of the vectors $\bfc_j$ can be found separately:
\begin{equation}\label{linsysc}
\bfD_{j} \bfc_j^* = \bfP_j - \bfB_j \bfa - \bfC_j \bfb^*,
\end{equation}
where (for coalescence) $\bfD_{j} = \bfV$ is a Vandermonde matrix \cite{nr} formed from the volume abscissas $v_n$. Other choices of moments have also been found to be numerically stable.  For example, another valid choice is
\begin{align}
\nonumber
& k=(i-1)/3 \quad i\in \{1,\ldots,2N\} \quad &\mbox{with} \quad &l=m=p=0 \\
\nonumber
& k=(2i-1)/3 \quad i\in \{1,\ldots,N\} \quad &\mbox{with}  \quad &l=1, \quad m=p=0\\
\nonumber
& k=(2i-1)/3 \quad i\in \{1,\ldots,N\} \quad &\mbox{with} \quad &m=1, \quad l=p=0\\
& k=(2i-1)/3 \quad i\in \{1,\ldots,N\} \quad &\mbox{with} \quad &p=1, \quad l=m=0.
\label{set3}
\end{align}
This choice can be found to give more accurate results for our test cases and still includes the surface area and the volume of the droplets, as well as their momentum.
Thus, it will be used for the computations in Section~\ref{results}.
We should note that for a given value of $N$, the simulation results found with the moment sets in Eqs.~(\ref{set3i}) and (\ref{set3}) were nearly identical.  The choice between these two systems was thus made based on ease of solution of the linear system.

In general, when moments involving the velocity are limited to first order, the matrices that must be inverted will be non-singular as long as the volume abscissas are distinct.
The numerical treatment of the singularities associated with Eq.~(\ref{linsys2}) has been discussed elsewhere \cite{marchisio05}.
The coalescence operator will normally force the $v_n$ to remain distinct if they have distinct velocities.
However, if due to initial conditions two or more of the volume abscissas are equal, it suffices to perturbate the values of $v_n$ enough to allow for the coefficient matrix in Eq.~(\ref{linsys2}) to be invertible.
We should also note that for cases dominated by coalescence (e.g., without evaporation) the volume abscissas grow rapidly, leading to matrices that are more and more ill-conditioned as the abscissas increase.
Thus, even though they are strictly non-singular, such matrices lead to severe numerical difficulties.
However ill-conditioning can be almost completely alleviated by using iterative improvements of the linear solver as described in Section 2.5 of Press et al.~\cite{nr}
after rescaling Eq.~(\ref{u3ns}).
The latter is done by defining positive scaling factors $v_{\rm s}$ and $u_{\rm s}$, and dividing both sides of Eq.~(\ref{u3ns}) by $v_{\rm s}^k u_{\rm s}^{l+m+p}$.
Note that the abscissas and unknown source terms are rescaled in a consistent manner: $v_n \rightarrow v_n /v_{\rm s}$, $\bfu_{n} \rightarrow \bfu_{n} /u_{\rm s}$, $\bfa \rightarrow \bfa$, $\bfb^* \rightarrow \bfb^*/v_{\rm s}$, and $\bfc_j^* \rightarrow \bfc_j^*/(v_{\rm s} u_{\rm s})$.
The evaporative flux constraints (Eqs.~(\ref{rca}), (\ref{rcv}) and (\ref{evap2})) can be rescaled in a similar manner by introducing a positive scaling factor $w_{\rm s}$ for the weights: $w_n \rightarrow w_n /w_{\rm s}$.
In this work, we use the following scaling factors: $v_{\rm s} = \max_n v_n$, $u_{\rm s} = \max_n |\bfu_n|$ and $w_{\rm s} = \sum_n w_n$.
We find that using the scaled variables and at most three iterative improvements of the linear solver are enough to completely eliminate round-off error in the solution to the DQMOM linear system.
Moreover, because round-off error leads to poor performance of the differential equation solver, the overall computational cost using the iterative improvements can be significantly reduced.

\section{Nozzle Test Problem}\label{test}

In order to validate the proposed DQMOM approach for strongly coalescing sprays, and to compare this method to both a reference Lagrangian solver solution as well as the solution obtained with the multi-fluid model, we need a well-suited test problem that is difficult enough to highlight the limitations of the methods under consideration.
For that purpose, we have chosen for the gas phase a 2D axisymmetrical conical decelerating nozzle,
designed in such a way that it admits, for one-way coupling spray dynamics a self-similar solution.
After presenting the details of this configuration, we will provide the set of DQMOM equations to be solved in this framework. We have selected six representative test cases, combining coalescence/no coalescence with evaporation/no evaporation, which are then presented.
Next we give an overview of
the Lagrangian solver that provides the reference solution for the various test cases. Because the problems under consideration can be difficult to solve numerically, we must be very careful as far as this reference solution in concerned and thus we provide the details of the Lagrangian numerical integration in the limit of one-way coupling with the gas phase.
Finally, before discussing the results in Section~\ref{results}, we review the fundamentals of the Eulerian multi-fluid model for the sake of self-consistency of the paper.

\subsection{Definition of configuration}

The chosen configuration is stationary 2D axisymmetrical in space and 1D in droplet size. 
It is described in detail, along with the Lagrangian solver, in \cite{lmv04}.
Hence, only its essential characteristics are given here.

A spray of pure heptane fuel is carried by a gaseous mixture of heptane and nitrogen into a conical diverging nozzle of axis ($0<z$).
At the entrance, $z=z_0$, 99\% of the mass of the fuel is in the liquid phase, whereas 1\% is in the gaseous mixture.
The temperature (400~K) as well as the composition of the gas mixture (mass fraction is 2.9\% for heptane and 97.1\% for nitrogen) is fixed during the entire calculation. The gas density is then $0.871733$~mg.cm${}^{-3}$.
The influence of the evaporation process on the gas characteristics is not taken into account in our one-way coupled calculation.
It is clear that the evaporation process is going to change the composition of the gas phase and then of the evaporation itself.
However, we do not attempt to achieve a fully coupled calculation, but only to compare two ways of evaluating the coupling of the dynamics, evaporation and coalescence of the droplets.
It has to be emphasized that it is not restrictive in the framework of this study, which is focused on the numerical validation of Eulerian solvers for the liquid phase under conditions of strong coalescence.

For the problem to be  one-dimensional in space, conditions for straight trajectories are used and are compatible with the assumption of an incompressible gas flow.
This leads to the following expression for the gaseous axial velocity $v_z$ and the reduced radial velocity $v_r/r$:
\begin{equation}\label{fluidv}
v_z=V(z)=\dfrac{z_0^2 V(z_0)}{z^2}, \quad
\frac{v_r}{r}=U(z)=\dfrac{V(z)}{z}=\dfrac{z_0^2V(z_0)}{z^3} \quad
\text{for $z\ge z_0$}
\end{equation}
where $z_0>0$ is the coordinate of the nozzle entrance and the axial velocity $V(z_0)$ at the entrance is fixed. The trajectories of the droplets are also assumed straight since their injection velocity is co-linear to the one of the gas. This assumption is only valid when no coalescence occurs. However, even in the presence of coalescence, it is valid in the neighborhood of the centerline.

\begin{figure}[t]
\noindent
\begin{center}
\epsfxsize=5in
\epsffile{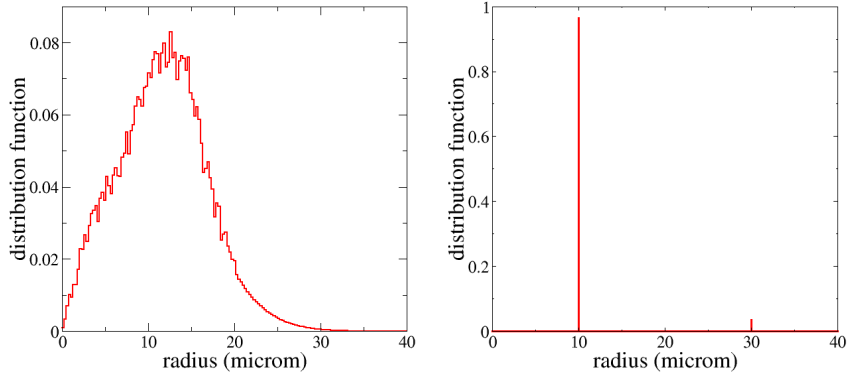}
\end{center}
\caption{Initial number density functions for droplet radius. Left: Monomodal distribution. Right: Bimodal distribution.}
\label{fig1}
\end{figure}

Let us finally consider two droplet distribution functions. The first one, called monomodal, is composed of droplets with radii between 0 and 35~$\mu$m, with a mean radius of 12~$\mu$m, a variance of 5~$\mu$m and a Sauter mean radius of 15.6~$\mu$m.
It is represented in Fig.~\ref{fig1} and is typical of the experimental conditions reported in \cite{lsgsm04}.
The droplets are constituted of liquid heptane, their initial velocity is the one of the gas, their initial temperature, fixed at the equilibrium temperature 325.4~K (corresponding to an infinite conductivity model), does not change along the trajectories.
The second distribution is called bimodal since it involves only two groups of radii, respectively, $10$ and $30$~microns with equal mass density.
This bimodal distribution function is typical of alumina particles in solid propergol rocket boosters \cite{hylkema99}, and is represented in Fig.~\ref{fig1}.

The initial injected mass density is taken as $m_0= 3.609$~mg.cm${}^{-3}$ so that the volume fraction occupied by the liquid phase is 0.57\%.
Because of the deceleration of the gas flow in the conical nozzle, droplets are going to decelerate, however at a rate depending on their size and inertia.
This will induce coalescence. The deceleration at the entrance of the nozzle is taken as $a(z_0)=-2 V(z_0)/z_0$; it is chosen large enough so that the velocity difference developed by the various sizes of droplets is important.
We have chosen rather large values, as well as strong deceleration, leading to extreme cases: $V(z_0)=5$~m.s${}^{-1}$, $z_0=10$~cm for the monomodal case and $V(z_0)=5$~m.s${}^{-1}$, $z_0=5$~cm for the bimodal case.
These values generate a very strong coupling between coalescence, evaporation and droplet dynamics.
These severe conditions as well as the two types of size distributions make the test cases under consideration very efficient tools for the numerical evaluation of the two Eulerian models.

\subsection{DQMOM model equations in nozzle configuration}

For the nozzle test case, Eqs.~(\ref{number}--\ref{momentum}) reduce to a set of ordinary differential equations (ODEs) defined on the interval $z \in [z_0, \infty)$ for the variables $w_n$, $v_n$, $\eta_n (z) = u_r / r$ and $\xi_n (z) = u_z$:
\begin{align}
2 w_n \eta_n + \partial_{z} \left( w_n \xi_n \right)
&= a_n , \label{numbera}   \\
2 w_n v_n \eta_n +  \partial_{z} \left( w_n v_n \xi_n \right)
&= b_n , \label{massa}  \\
3 w_n v_n \eta_n^2  + \partial_{z} \left( w_n v_n \eta_n \xi_n \right)
&= c_{rn}/r , \label{momentuma1} \\
2 w_n v_n \eta_n \xi_n  + \partial_{z} \left( w_n v_n \xi_n^2 \right)
&= c_{zn} , \label{momentuma2}
\end{align}
where $u_z = \xi(z)$ and $u_r = r \eta(z)$ are the axial and radial components of the spray velocity, respectively. The corresponding fluid velocities are given in Eq.~(\ref{fluidv}).
The terms on the right-hand side of Eqs.~(\ref{massa}--\ref{momentuma2}) are given by
\begin{align}
b_n &= b^*_n + w_n R_v (v_n) , \label{bn}  \\
c_{rn}/r &= c^*_{rn} + w_n \eta_n R_v(v_n)
+ w_n v_n F_r(v_n,\eta_n)/r, \label{cn1} \\
c_{zn} &= c^*_{zn} + w_n \xi_n R_v(v_n)
+ w_n v_n F_z(v_n,\xi_n), \label{cn2}
\end{align}
where the drag model is
\begin{align}
F_r(v,\eta)/r &= \alpha \left( \frac{4 \pi}{3v} \right)^{2/3} (U - \eta), \label{dn1} \\
F_z(v,\xi) &= \alpha \left( \frac{4 \pi}{3v} \right)^{2/3} (V - \xi) \label{dn2}
\end{align}
with $\alpha = 1.566 \times 10^{-7}$~m$^2$.s${}^{-1}$.

From the form of the governing equations, it is straightforward to show that if $\eta_n = \xi_n/z$ at $z=z_0$, then this relation will hold for all $z$ and the droplet trajectories are straight lines.
The system of DQMOM model equations can thus be reduced to three nonlinear ODEs for $w^*_n = w_n (z/z_0)^2$, $v_n$, and $\xi_n$ by eliminating Eq.~(\ref{momentuma1}):
\begin{align}
\partial_{z} \left( w^*_n \xi_n \right)
&= a_n , \label{numbere}   \\
\partial_{z} \left( w^*_n v_n \xi_n \right)
&= b^*_n + w^*_n R_v (v_n) \label{masse}  \\
\intertext{and}
\partial_{z} \left( w^*_n v_n \xi_n^2 \right)
&= c^*_{zn} + w^*_n \xi_n R_v (v_n)
+
\alpha w^*_n v_n \left( \frac{4 \pi}{3v_n} \right)^{2/3} (V - \xi_n).
\label{momentum2e}
\end{align}
The terms on the left-hand side represent changes in the weights and abscissas due to transport.  The terms on the right-hand side represent, respectively, the changes due to coalescence, evaporation, and drag. The coalescence terms are found by solving
\begin{multline}
 (1-k) \sum_{n=1}^{N} v_n^{k} \xi_{n}^{m} a_n
+ (k - m) \sum_{n=1}^{N} v_n^{k-1} \xi_{n}^{m} b_n^*
+ m \sum_{n=1}^{N} v_n^{k-1} \xi_{n}^{m-1} c_{zn}^*  \\
=
\frac{1}{2} \left( \frac{z_0}{z} \right)^2
\sum_{n=1}^{N} \sum_{q=1}^{N}
w^*_n w^*_q B ( |\xi_n - \xi_q |, v_n , v_q ) \\
\left[ ( v_n + v_q )^k
\left( \frac{v_n \xi_{n} + v_q \xi_{q}}{v_n + v_q} \right)^{m}
- v_{n}^k \xi_{n}^{m}
- v_{q}^k \xi_{q}^{m}
\right] . \label{coalne}
\end{multline}
Note the presence of the scaling factor $(z_0/z)^2$ in the coalescence rate. As discussed in Section~\ref{decomp}, we will use moments given in Eq.~(\ref{set3}) that decouple Eq.~(\ref{coalne}) into two smaller systems.

\subsection{Test cases}

For evaporation, we will consider three cases described below: (\emph{i}) no evaporation ($R_v=0$), (\emph{ii}) linear evaporation ($R_v \propto v$), and (\emph{iii}) non-linear evaporation ($R_v \propto v^{1/3}$).
For each case, we will consider two sub-cases: without coalescence ($E_{\rm coal}=0$) and with coalescence ($E_{\rm coal}=1$).
The two evaporation laws correspond to the two cases described in Section~\ref{ef}, for which droplets disappear either in infinite time (\emph{ii}), thus leading to a evaporative flux at zero size, or in finite time (\emph{iii}) for which the evaporative flux depends on the structure of the number density function in size phase space.

\subsubsection*{No evaporation}

For the special case of no evaporation and no drag, the right-hand sides of Eqs.~(\ref{massa}--\ref{momentuma2}) are null.
This special case has an analytical solution with $v_n$, $w^*_n$, and $\xi_n$ constant.  In the opposite limit of no evaporation and infinite drag, $\xi_n =  V$ and $w^*_n \propto (z/z_0)^2$.

For non-evaporating droplets, $R_v =0$. In the absence of coalescence, $a_n=b^*_n=c^*_{zn}=0$.
The DQMOM model reduces to $v_n$ and $w^*_n$ constant, and
\begin{equation}\label{xin}
\xi_n \partial_{z} \xi_n =
 \alpha \left( \frac{4 \pi}{3v_n} \right)^{2/3} (V - \xi_n).
\end{equation}
This result is consistent with our earlier remark concerning the cases of zero and infinite drag.  Finally, we should note that even with coalescence the momentum is conserved ($k=m=1$) so that $\sum c^*_{zn}=0$.
Thus, we can expect $w^*_n \xi_n$ to be approximately constant for all values of drag.  For this case we expect excellent agreement between DQMOM and the Lagrangian solver in the absence of coalescence since the corresponding transport equations are identical (i.e., each DQMOM abscissa behaves like a Lagrangian particle).
On the other hand, with coalescence the droplets grow very large and we expect differences due to how the coalescence term is treated in each method.
This test case will, however, be very difficult for the multi-fluid model, since it was especially designed to tackle the problem of evaporation.
In the presence of strong growth of droplet size, the number of sections that must be used in order to reproduce the physics with the multi-fluid model will also dramatically increase.
Consequently, this test case will allow us to both test the capability of the DQMOM to capture the coupling of dynamics and coalescence at low cost in comparison to the Lagrangian solution, and to see if the multi-fluid model can provide good results, even if it  is not competitive in terms of computational efficiency.

\subsubsection*{Linear evaporation}

For evaporating droplets with linear evaporation, we take
\begin{equation}
R_v (v_n) = - E_v v_n,
\end{equation}
with $E_v= 7.1262$~s${}^{-1}$ for the monomodal case and $E_v= 14.2524$~s${}^{-1}$ for the bimodal case.
For this case, the evaporative flux $\psi$ is zero. The coalescence terms are again found by solving Eq.~(\ref{coalne}).
In the absence of coalescence, we have $a_n=b^*_n=c^*_{zn}=0$ and the DQMOM model reduces to $w^*_n \xi_n$ constant, Eq.~(\ref{xin}), and
\begin{equation}\label{vn}
\xi_n \partial_{z} v_n = R_v (v_n) .
\end{equation}
Thus the volume $v_n$ and velocity $\xi_n$ are coupled through evaporation and drag, but are independent of $w^*_n$ in the absence of coalescence.
For this case we again expect excellent agreement between DQMOM and the Lagrangian solver in the absence of coalescence since the corresponding transport equations are identical.
On the other hand, with coalescence there is a competition between growth and evaporation leading to smaller droplets than in the non-evaporating case.
This is a very interesting test case, since it will allow us to compare both methods in an evaporative configuration, but without getting into the difficulty of modeling the droplet disappearance with the DQMOM approach.

\subsubsection*{Non-linear evaporation}

With non-linear evaporation we will use
\begin{equation}
R_v (v_n) = - \frac{E_{s}}{2} \left( \frac{3 v_n}{4 \pi} \right)^{1/3} \label{nle}
\end{equation}
with $E_{s}=1.99 \times 10^{-7}$~m$^2$/s.  For this case the evaporative flux $\psi$ will generally be non-zero, and is found using the method described in Section~\ref{ef} with $w^*_n$ in place of $w_n$. However, we will also compare predictions for the bimodal initial distribution found by setting $\psi=0$. As for the previous cases, we will investigate the effect of the flux model with and without coalescence.  From a practical standpoint, the behavior of DQMOM with non-linear evaporation is of great interest and it is a configuration with which the comparison of both Eulerian models will be of practical relevance.

\subsection{Reference Lagrangian solution}

Euler-Lagrange numerical methods are commonly used for the calculation of polydisperse sprays in various application fields (see for example \cite{orourke81,hylkema98,miller99,miller00,ruger00,dupays00} and the references therein).
In this kind of approach, the gas phase is generally computed using a deterministic
Eulerian solver, while the dispersed phase is treated 
in a Lagrangian way.
The influence of the droplets on the gas flow is taken into account by the presence of source terms 
in the system of gas conservation equations. 
Two Lagrangian methods can be used as far as the dispersed phase is concerned depending on the level at which the physical processes are modeled. The first one is a 
Discrete Particle Simulation, which describes the evolution of numerical particles, each one representing one or several droplets. The physical processes such as transport, evaporation, collisions are then described by Liouville equations and the Eulerian fields usually recovered through ensemble averages. However, in the present study, we have preferred  the Williams governing equation 
and thus a statistical description of the coalescence process. We then coherently use a 
Direct Simulation Monte-Carlo method (DSMC), the second kind of approach.
It can be seen as the uncoupling, over a small time step, of the droplet transport in phase space (dynamics and evaporation), described by a particle method, and the collisions described by a Monte-Carlo method.

 A  complete exposition on the derivation and implementation of this method is outside the scope of this paper. We refer the reader, for example, to \cite{amsden89,hylkema98,hylkema99} for more details.   
 Here, for the sake of completeness, we present only the main features of the numerical method that we used in order to provide a ``reference numerical solution'' for the test cases.

\subsubsection*{Lagrangian solver}


The Lagrangian solver can be roughly interpreted as a stochastic representation of the kinetic equation (\ref{weq}). In other words, in the limit of a sufficiently large number of stochastic particles and a sufficiently fine computational grid (at least in the case of one-way coupling), the statistical estimates for the moments found from the particles should converge to those computed from the Eq.~(\ref{weq}).
In the Lagrangian solver, at each time step $k$, the droplet distribution function $f(t^k)$ is approximated by a finite
 weighted sum of Dirac masses, $\tilde{f}(t^k)$, which reads
\begin{equation}
\tilde{f}(t^k) = \sum_{i=1}^{N^k} n_i^k  \, \delta_{z_i^k, u_i^k, v_i^k}.
\end{equation}
Each weighted Dirac mass is generally called a ``parcel" and
 can be physically interpreted as an aggregated number of droplets (the weight $n_i^k$),  located around the same  point, $z_i^k$,  with about the same  velocity,
  $u_i^k$ and  about the same volume, $v_i^k$.  $N^k$ denotes the total number of parcels in the computational domain at time $t^k$.
In all our calculations, the weights $n_i^k$ were chosen in such a way that each parcel represents the same volume of liquid ($n_i^k v_i^k = Const$).


Each time step of the particle method is divided in two stages. The first is devoted to discretization of the left-hand side of the kinetic equation (\ref{weq}), modeling the motion and evaporation of the droplets.
 In our code, the new position, velocity
 and volume of each  parcel are calculated according to the following numerical scheme:
\begin{align}
\label{eq:syst1}
u_i^{k+1} & =  u_i^{k} \exp({-\Delta t/ \tau_i^k})
+ V(z_i^{k})  \bigl(1 - \exp({-\Delta t/ \tau_i^k})\bigr)  \\
\label{eq:syst2}
\int_{v_i^k}^{v_i^{k+1}} &dv/R_v(v) =   \Delta t\\
\label{eq:syst3}
z_i^{k+1} & =   z_i^k + \Delta t \, u_i^{k+1} =   z_i^k + \Delta t \, V(z_i^{k})
+ \Delta t\,(u_i^{k}  -   V(z_i^{k})) \exp({-\Delta t/ \tau_i^k})
\end{align}
where $V$ denotes the axial gas velocity,
$z_i^k$ ($u_i^k$) corresponds to the axial coordinate of the position (velocity) of the parcel $i$ at time $t^k$,
$R_v$ is the evaporation rate (independent of $t$ and $x$ since the gas composition and temperature are assumed constant in the domain).
Eq.~(\ref{eq:syst2}) is resolved analytically and depends on the chosen evaporation model.
For linear evaporation, it can be written as
\begin{equation}
v_i^{k+1} = v_i^{k} \exp(-E_v \Delta t)
\end{equation}
and for non-linear evaporation, it is written as
\begin{equation}
s_i^{k+1} = s_i^{k} -E_s \Delta t
\end{equation}
where $s_i^k$ is the parcel surface area.
The parcel relaxation time $\tau_i^k$ is defined as
\begin{equation}
\tau_i^k = \frac{ 2 \rho (r_i^k)^2 }{9 \mu_g},
\end{equation}
with $r_i^k$ being the parcel radius, $\rho$ the liquid density and
 $\mu_g$ the gas viscosity.

In Eqs.~(\ref{eq:syst1}--\ref{eq:syst3}), the  parcel radial coordinate is not calculated
because the trajectories of the parcels are straight lines.
Besides, as  mentioned above,  the influence of the droplets on the gas flow is not taken into account.
Hence, Eq.~(\ref{fluidv}) is used to calculate the gas velocity, $V(z_i^{k})$,  at the parcel location.


The second stage of a time step is devoted to the discretization of the collision operator.
Several Monte-Carlo algorithms have been proposed in the literature for the treatment of droplet collisions \cite{orourke81,hylkema98,ruger00,hylkema99,schmidt00}.
They are all inspired by the methods used in molecular gas dynamics \cite{bird94} and, in particular,  they  suppose that the computational domain is divided into cells, or control volumes, which are small enough to consider that, within them, the droplet distribution function is almost uniform.

 The algorithm used in our reference Lagrangian solver is close to the one proposed by O'Rourke.
 It consists of the following 3 steps (see also \cite{hylkema98} for more details):
\begin{itemize}
\item[1.] For each computational cell $C_J$, containing $N_J$ parcels,  we
 choose randomly,   with  a uniform distribution law,
 $N_J/2$  pairs of parcels, $(N_J-1)/2$ if $N_J$ is odd.

\item[2.]
For each pair $p$, let $p_1$ and $p_2$ denote the two corresponding parcels with the convention $n_1 \ge n_2$, where $n_1$ and $n_2$ denote the parcel numerical weights.
Then for each  pair $p$ of the cell $C_J$, we choose randomly  an integer $\nu_p$, according to the Poisson distribution law:
$$ P(\nu) = \frac{\lambda_{12}}{\nu !} \exp{(-\lambda_{12}) },$$
with
$$ \lambda_{12} =
\pi \frac{n_1 (N_J - 1) \Delta t}{{\rm vol}(C_J)} (R_1 + R_2)^2 |u_1 - u_2| $$
with ${\rm vol}(C_J)$ being the volume of cell $C_J$, which is proportional to $(z_J/z_0)^2$ for the nozzle test case, and $R_1$, $R_2$ being the radii of the parcels $p_1$, $p_2$.
The coefficient $\lambda_{12}$ represents the mean number of collision, during $(N_j - 1)$ time steps, between a given droplet of parcel $p_2$ and any droplet of parcel $p_1$. Note that a given pair of parcels is chosen, on average, every $(N_j - 1)$ time steps.

\item[3.] If $\nu_p = 0$, no collisions occur during this time step between the parcels $p_1$ and $ p_2$.
Otherwise, if $\nu_p > 0$, parcel $p_1$ undergoes $\nu_p$ coalescence with parcel $p_2$ and the outcome of a collision is treated as follows.
First the weight $n_1$ of the parcel  $p_1$ is replaced by $n'_1 = n_1 - \nu_p n_2$ and its other characteristics are left unchanged.
If $n'_1 \le 0$, parcel $p_1$ is removed from the calculation.
Secondly, the velocity $u_2$ and the volume $v_2$ of parcel $p_2$ are replaced by
$$ v_2' = v_2 + \nu_p  v_1,   \quad \quad  u'_2 =  \frac{ v_2 u_2 +  \nu_p v_1 u_1}{ v_2 + \nu_p  v_1}, $$
and its weight, $n_2$,  is left unchanged.
\end{itemize}

Let us mention that,  for each time step and each control volume $C_J$,  the computational cost of this algorithm scales like $O(N_J)$.
This is a great advantage compares to the O'Rourke method, which scales like
 $O(N_J^2)$.
 Another algorithm, with the same scaling features, has been introduced by Schmidt and Rutland in \cite{schmidt00}.


To obtain good accuracy, the time step, $\Delta t$, must be chosen small enough to ensure  that the number of collisions between two given parcels, $p_1$ and $p_2$, is such that for almost every time: $\nu_p n_2 \le n_1$.
The average value of $\nu_p$ being $\lambda_{12}$, this constraint is equivalent to the condition

\begin{equation}\label{limdt}
\frac{n_2 N_J \Delta t}{{\rm vol}(C_J)} \pi (R_1 + R_2)^2 |u_1 - u_2|  \ll  1.
\end{equation}

For the nozzle test case described above, this constraint reveals to be less restrictive
 than the ``CFL" condition
\begin{equation}\label{LAcfl}
\forall i=1, .... N, \quad   \frac{ |u_i | \Delta t } {\Delta z} \ll  1,
\end{equation}
with $\Delta z$ being the mesh size.
This  condition is necessary to compute accurately the droplet movement and in particular to avoid that a parcel goes through several control volumes during the same time step.
This is essential in order to have a good representation of the droplet distribution function  in each mesh cell.

\subsubsection*{Reference solution}

The Lagrangian solver just described is used to provide reference solutions in stationary cases with and without coalescence.
In order to obtain a converged solution, particular attention must be devoted to the choice of the number of parcels, the size of the cells, and the time step.

For cases without coalescence, the computational cells are only used to have spatially averaged quantities to compare with Eulerian results.
Moreover, the stationary aspect of the problem allows averaging in time in order to obtain smooth solutions.
For these reasons, the conditions on the number of parcels and on the size of the computational cells are not very restrictive in the absence of coalescence.
The time step is only limited by the CFL-like condition (\ref{LAcfl}) needed for the convergence, with a low value.
This last condition is the most restrictive since the scheme used for the transport of the particles is first order.
For our test cases, the time step must be $10^{-6}$~s or smaller.

\bigskip

\begin{table}[tbh]
\smallskip
\centering
\begin{tabular}{|c|c|c|c|}
\hline
distribution & evaporation & No. of parcels & No. of parcels inj./s \\
\hline
monomodal & linear & 41,560 & 100,000\\
monomodal & nonlinear & 20,440 & 1,000,000\\
bimodal & nonlinear & 6,320 & 200,000\\
\hline
\end{tabular}
\smallskip
\caption{Number of parcels for the Lagrangian simulations for the cases without coalescence.}\label{tableLass}
\end{table}

\medskip

\begin{table}[tbh]
\smallskip
\centering
\begin{tabular}{|c|c|c|c|c|}
\hline
distribution & evaporation & No. of parcels & No. of parcels inj./s & min. No. of          \\
             &             &                &                       & parcels /cell        \\
\hline
monomodal & linear & 160,000 & 200,000 & 40 \\
bimodal & linear & 126,000 & 560,000 & 50 \\
monomodal & nonlinear & 35,000 & 1,300,000 & 260 \\
monomodal & no & 44,200 & 300,000 & 65 \\
\hline
\end{tabular}
\smallskip
\caption{Number of parcels for the Lagrangian simulations for the cases with coalescence.}\label{tableLacol}
\end{table}

\medskip

For cases with coalescence, there are additional restrictions.
First, the algorithm used for coalescence assumes that the droplet distribution function of the spray is nearly uniform over each computational cell.
However, in the region with high gradients of the gas velocity, that is to say at the entrance of the nozzle, this distribution can change quickly and the size of the cells must be small enough to avoid numerical errors.
Moreover, in order to properly describe the coalescence phenomenon in each cell with the stochastic algorithm, a sufficient number of parcels must be present in each cell, typically on the order of 50, with a minimum of 20 \cite{alexander97}.
The smaller are the cells, the larger must be the number of parcels in the computational domain.
The required size of the cells is evaluated for the case where the size distribution changes the most rapidly (the case without evaporation).
We then employ a non-uniform space discretization with 130 cells, with smaller cells near the entrance of the nozzle defined using a uniform discretization for the variable $z^{3/10}$.
The number of parcels injected per second is given in Tables~\ref{tableLass} and \ref{tableLacol}.

\subsection{Eulerian multi-fluid solver}

Eulerian multi-fluid methods were developed as an alternative to Lagrangian methods for the simulation of polydisperse evaporating sprays.
A complete derivation of such methods from the kinetic model is performed in \cite{laurent01} for dilute sprays and in \cite{lmv04} for sprays with coalescence.
The principle of the method is quite different than the one used in DQMOM.
Indeed, it can be considered as a finite-volume discretization in the droplet size phase space for moments of order $0$ and $1$ of the velocity distribution conditioned on size.

In laminar flows, it can be proven rigorously that it is sufficient to work with only these two moments as long as the velocity distribution conditioned on droplet size is mono-kinetic \cite{dufour_th,dufour05} (i.e.\ all droplets with the same volume have identical velocities so that the size-conditioned velocity variance is null.)  By construction, the nozzle test problem will be mono-kinetic for non-coalescing droplets. However, with coalescence the size-conditioned velocity variance can be nonzero. Comparisons between the Lagrangian and Eulerian results in the presence of coalescence will therefore allow us to quantify the magnitude of the error caused by invoking the mono-kinetic assumption in the Eulerian models. (Recall that the choice of moments used in the DQMOM linear system is equivalent to the mono-kinetic assumption in the multi-fluid model.)
In this section, we provide only the main points of the derivation of the multi-fluid model, as well as the underlying assumptions that are implied, and the resulting system of equations that will be solved.

The first step consists of writing equations for the two moments in velocity.
This leads to the closed semi-kinetic model if the following assumption is made concerning the structure of $f$:
$f (v, \bfu; \bfx, t) = n(v; \bfx, t)\delta(\bfu-\bar\bfu(v; \bfx, t))$.
In other words, the droplet velocity conditioned on the size is assumed to be Dirac delta function.
In the case of a coalescing spray, the compatibility of such a condition with droplet
coalescence is far from obvious; however, the semi-kinetic system of conservation equations can be obtained
by using an asymptotic limit as presented in \cite{lmv04}.

The second step consists of discretizing $n(v)$ in sections $[v^{(j-1)},v^{(j)})$ and
in integrating the semi-kinetic model over each section.
This leads to a multi-fluid model (by using a presumed distribution $\kappa^{(j)}(v)$ in each section), thereby yielding a conservation equation on the moment associated with the mass density
$$
n(v; \bfx, t) = m^{(j)}(\bfx, t)\kappa^{(j)}(v) \quad \mbox{where} \quad
\int_{v^{(j-1)}}^{v^{(j)}}\rho v\kappa^{(j)}(v)dv = 1.
$$
In addition, only the averaged velocity is considered in each section, i.e.\
$\bar\bfu(v; \bfx, t) = \bar\bfu^{(j)}(\bfx, t)$ if $v^{(j-1)} \le v < v^{(j)}$.
The resulting system can be found in \cite{lmv04}.
It can be rewritten and simplified in the stationary, self-similar, 2D axisymmetrical configuration we are considering.

The resulting set of equations is
\begin{equation}\label{eq:mj}
2 m^{(j)} \frac{{u_z}^{(j)}}{z} +\partial_z (m^{(j)} {u_z}^{(j)})
=-(E_1^{(j)}+E_2^{(j)})m^{(j)}+E_1^{(j+1)}m^{(j+1)} + C_{\rm m}^{(j)}
\end{equation}
\begin{multline}\label{eq:mjuj}
2 m^{(j)} \left(\frac{{u_z}^{(j)}}{z}\right)^2 +\partial_z (m^{(j)} ({u_z}^{(j)})^2)
= \\
-(E_1^{(j)}+E_2^{(j)})m^{(j)}{u_z}^{(j)}
+E_1^{(j+1)}m^{(j+1)}{u_z}^{(j+1)}
+ m^{(j)} F_z^{(j)}+C_{\rm muz}^{(j)}
\end{multline}
where ${u_z}^{(j)}$ is the axial velocity, which only depends on $z$, and $r{u_z}^{(j)}/z$
is the radial velocity, since the trajectories are straight lines.
Moreover, $E_1^{(j)}$ and $E_2^{(j)}$ are the
``classical'' pre-calculated vaporization coefficients
\cite{greenberg93,laurent01}:
\begin{equation*}
E_1^{(j)} = -\rho\, v^{(j-1)}\, R_v(v^{(j-1)})\,
\kappa^{(j)}(v^{(j-1)})\quad \mbox{and} \quad
E_2^{(j)} = -\int_{v^{(j-1)}}^{v^{(j)}}
\rho\,R_v(v)\,\kappa^{(j)}(v)\,{\rm d}v,
\end{equation*}
and $F_z^{(j)} (v_u^{(j)},{u_z}^{(j)})$ is the axial component of the ``classical'' pre-calculated
drag force \cite{greenberg93,laurent01}:
\begin{equation*}
F_z^{(j)}=
\int_{v^{(j-1)}}^{v^{(j)}}\rho\,v\,
F_z(v,{u_z}^{(j)})\,\kappa^{(j)}(v)\,{\rm d}v
\quad \mbox{where} \quad
v_u^{(j)} = \left[
\frac{\int_{v^{(j-1)}}^{v^{(j)}}v\,\kappa^{(j)}(v)\,{\rm d}v}
{\int_{v^{(j-1)}}^{v^{(j)}}v^{1/3}\,\kappa^{(j)}(v)\,{\rm d}v}
\right]^{3/2}.
\end{equation*}

The source terms associated with coalescence phenomenon in the mass and momentum equation, respectively, of the $j$th section read
\begin{equation}\label{eq:cmj}
C_{\rm m}^{(j)} =
-m^{(j)} \sum\limits_{k=1}^{N} m^{(k)}V_{jk} Q_{jk}
+\sum\limits_{i=1}^{I^{(j)}} m^{(o^\diamond_{ji})} m^{(o^*_{ji})}
V_{o^\diamond_{ji}o^*_{ji}} (Q_{ji}^{\diamond}+Q_{ji}^{*}),
\end{equation}
\begin{multline}\label{eq:cmuj}
C_{\rm muz}^{(j)} =
-m^{(j)} {u_z}^{(j)}\sum\limits_{k=1}^{N} m^{(k)} V_{jk} Q_{jk} \\
+\sum\limits_{i=1}^{I^{(j)}}m^{(o^\diamond_{ji})} m^{(o^*_{ji})}V_{o^\diamond_{ji}o^*_{ji}}\left({u_z}^{(o^\diamond_{ji})}Q_{ji}^{\diamond}+{u_z}^{(o^*_{ji})}Q_{ji}^{*}\right),
\end{multline}
where $V_{jk}= | {u_z}^{(j)}-{u_z}^{(k)}|$ and the collision integrals
$Q_{jk}$, $Q_{ji}^{\diamond}$ and $Q_{ji}^{*}$ do not depend on $z$.
The disappearance integrals  ${Q}_{jk}$
are evaluated on rectangular domains
$L_{jk} = [v^{(j-1)},v^{(j)}]\times[v^{(k-1)},v^{(k)}]$,
whereas the appearance integrals,
 $Q_{ji}^{\diamond}$ and $Q_{ji}^{*}$,
 are evaluated on the diagonal strips
$D^{\diamond*}_j = \{ (v^\diamond,v^*), v^{(j-1)}\le v^\diamond+v^*\le v^{(j)}\}/\cup_{k=1}^N L_{kk}$,
which are symmetric strips with respect to the axis
$v^\diamond = v^*$.
These strips $D^{\diamond*}_j$ are divided into domains, denoted by
$X_{ji}$ and the symmetric one, $X_{ji}^{{\rm sym}}$, where the velocity
of the partners is constant. The domains
$X_{ji}$ and $X_{ji}^{{\rm sym}}$ are the intersection of
$D^{\diamond*}_j$ with  $L_{kl}$, $k>l$, and $L_{kl}$, $k<l$, respectively;
their index is denoted $i \in [1,I^{(j)}]$ and
we define two pointers that indicate the collision
partners for coalescence, at fixed $i$:  $o^\diamond_{ji} = k$
and $o^*_{ji}=l$.

The coefficients used in the model, either
for the vaporization process or the drag force
$E_1^{(j)}$, $E_2^{(j)}$ and $F_z^{(j)}$, $j=[1,N]$ in Eqs.~(\ref{eq:mj}--\ref{eq:mjuj}),
or for the  coalescence:
$Q_{jk}^{}$, $j=[1,N],k=[1,N], k\ne j$, $Q_{ji}^{\diamond}$, $Q_{ji}^{*}$, $j=[2,N],i=[1,I^{(j)}]$ in Eqs.~(\ref{eq:cmj}--\ref{eq:cmuj}) can be pre-evaluated from the choice
of $\kappa^{(j)}$ in each section. The algorithms for the evaluation
of this coefficients are provided in \cite{lmv04}.
The distribution function is chosen
constant as a function of the radius in the sections 1 to $N$
and exponentially decreasing as a function of the surface in the last section, as done in \cite{lmv04}.

Because only the one-way coupled equations are solved and since the structure of the gas  velocity field is prescribed and stationary, we only have to solve the 1D ordinary differential Eqs.~(\ref{eq:mj},\ref{eq:mjuj}) for  each section.
The problem is then reduced to the integration of a stiff initial value problem from the inlet where the droplets are injected until the point where 99.9\% of the mass has evaporated.
The integration is performed using LSODE for stiff ordinary differential equations from the ODEPACK library \cite{hindmarsh83}.
It is based on BDF methods \cite{hairer91} (Backward Differentiation Formulae) where the space step is evaluated at each iteration, given relative and absolute error tolerances \cite{hindmarsh83}.
The relative tolerance, for the solutions presented in the following, are taken to be $10^{-4}$, 
and the absolute tolerance are related to the initial amount of mass in the various sections, since it can vary of several orders of magnitude.
Repeated calculations with smaller tolerances have proved to provide essentially the same solutions.

\section{Results and Discussion}\label{results}

Simulations for the cases presented in the previous section were carried out with the Lagrangian method, the multi-fluid model, and DQMOM.
Except for the cases without evaporation for which the multi-fluid method is not well suited (it requires a large number of sections and is only presented for comparison purposes),  the Eulerian methods were solved using a initial-value solver for ODEs and
required very short computational times (i.e., CPU secs) on a desktop computer. It is interesting to note that in the case without evaporation, the small computational cost still holds for the DQMOM approach.

In contrast, the time- and $z$-dependent Lagrangian simulations required several CPU hrs for each case.
Because the DQMOM and multi-fluid results do not depend on time, it is not appropriate to compare the computing times directly.
Nevertheless, it will generally be the case that using Eulerian methods will result in a substantial reduction
in the computing time for solving the spray equation. Such a statement was studied in details in \cite{lmv04} for unsteady
calculations and the conclusions drawn from that paper are applicable to the two Eulerian methods presented here.
Thus, the principal open question is whether or not the DQMOM results are of comparable
accuracy to the multi-fluid model and to the more costly Lagrangian  simulations.
We will compare predictions for selected statistics from the three solution methods in order to answer this question.
For the monomodal distribution and DQMOM resolution, several initial conditions will be used in the following an are presented in Table~\ref{table1}.

\begin{table}[tbh]
\smallskip
\centering
\begin{tabular}{|c|cc|cc|cc|cc|}
\hline
\multicolumn{9}{|c|}{Monomodal distribution}\\ \hline
    &\multicolumn{2}{|c|}{Vol. moments, $N$=4}
    &\multicolumn{2}{|c|}{Rad. moments, $N$=4}
    &\multicolumn{2}{|c|}{Rad. moments, $N$=6}
    &\multicolumn{2}{|c|}{Rad. moments, $N$=8}\\
$n$ & $w_n/N_0$   & $r_n$     & $w_n/N_0$ & $r_n$    & $w_n/N_0$  & $r_n$      & $w_n/N_0$     & $r_n$      \\ \hline
 1  & 0.7323      & 9.9955    & 0.1845   &  4.4079 & 8.5573E-2 & 3.3423 & 4.6445E-2 & 2.8465\\
 2  & 0.2545      & 18.5282   & 0.5397   & 11.0409  & 0.2779   & 7.5262 & 0.1488 & 5.5373\\
 3  & 1.288E-2    & 27.5630   & 0.2635   & 18.2840  & 5.5339E-2 & 12.9743 & 0.3089 & 9.6916\\
 4  & 2.279E-4    & 36.0142   & 1.212E-2 & 28.3910 & 4.9778E-3 & 18.8823 & 0.3438 & 14.2697\\
 5  &  & &  & & 3.1137E-4 & 26.3693 & 0.12931 & 19.2986\\
 6  &  & &  & & 1.6671E-5 & 34.7171 & 2.0905E-2 & 25.2866\\
 7  &  & &  & & & & 1.6982E-3 & 31.5808\\
 8  &  & &  & & & & 6.5627E-5 & 37.5149\\
 \hline
\end{tabular}
\smallskip
\caption{Initial conditions for weights and abscissas found using QMOM.}\label{table1}
\end{table}

\bigskip

The representative moments used to compare the three solution methods are the number density $m_0$, the mass density $m_1$, the average axial velocity difference between droplets and the gas phase $u_d$, and the Sauter mean radius $r_{3 2}$.
They are defined by
\[
m_0 = \int f(v, \bfu ) \sd v \sd \bfu ,
\qquad
m_1 = \int \rho v f(v, \bfu ) \sd v \sd \bfu,
\]
\[
u_d
= \dfrac{\int \rho v u_z f(v, \bfu ) \sd v \sd \bfu}{m_1} - V ,
\qquad
r_{3 2}
= \left( 4 \pi / 3 \right)^{1/3}
\dfrac{\int v f(v, \bfu ) \sd v \sd \bfu}{ \int v^{2/3} f(v, \bfu ) \sd v \sd \bfu} .
\]
With the DQMOM approach, these quantities are written
\[
m_0
= \sum_{n=1}^{N} w_n ,
\qquad
m_1
= \sum_{n=1}^{N} w_n \rho v_n ,
\]
\[
u_d
= \sum_{n=1}^{N} w_n v_n ( \xi_n - V) / m_1 ,
\qquad
r_{3 2}
= \left( 4 \pi / 3 \right)^{1/3}
\frac{\sum_{n=1}^{N} w_n v_n }{\sum_{n=1}^{N} w_n v_n^{2/3} }.
\]
And with the multi-fluid method, they are
\[
m_0
= \sum_{j=1}^{N} m^{(j)}\int_{v^{(j-1)}}^{v^{(j)}} \kappa^{(j)}(v) dv  ,
\qquad
m_1
= \sum_{j=1}^{N} m^{(j)},
\]
\[
u_d
= \sum_{j=1}^{N} m^{(j)}({u_z}^{(j)}-V) / m_1 ,
\qquad
r_{3 2}
= \left( 4 \pi / 3 \right)^{1/3}
\frac{\sum_{j=1}^{N}m^{(j)} \int_{v^{(j-1)}}^{v^{(j)}} v\kappa^{(j)}(v) dv  }
{\sum_{j=1}^{N} m^{(j)} \int_{v^{(j-1)}}^{v^{(j)}}v^{2/3} \kappa^{(j)}(v) dv }.
\]
Note that in practical applications, the mass density is a key quantity because it represents the total mass of liquid contained in the droplets.  In the nozzle test case, the rate of coalescence is strongly dependent on the velocity difference between droplets, which we find to be strongly correlated with the average axial velocity difference.  Indeed, if $u_d$ is not accurately captured, then we find that the predictions for all moments will degrade accordingly. In addition to the moments, we will also compare the mean droplet velocity conditioned on the radius $\langle u_z | r \rangle$ at selected downstream locations, as well as the mass distribution function ($\rho vf$). For the DQMOM, the scaled weights will be used to represent the mass distribution function. Obviously, since the sum of the weights equals the area under the mass distribution function, the absolute value of the heights of the scaled weights is arbitrary. Nevertheless, the relative heights and the locations provide insight into how well the quadrature points represent the distribution function.

We should note that for the monomodal cases without coalescence, the results with no evaporation were essentially identical for all three solution methods. The results presented below for the monomodal case with linear evaporation are representative of the quality of the predictions for all cases without coalescence and no evaporation. Likewise, for the bimodal case without coalescence and with linear or no evaporation, DQMOM and the Lagrangian method were essentially identical.  The multi-fluid method also yielded very good results for these cases if the number of sections was chosen large enough to mitigate
the numerical diffusion 
in the size phase space associated to the description of evaporation
that leads to broadening of the peaks.  Nonetheless, because none of these cases revealed any unanticipated problems for any of the simulation methods, we will not discuss them further.
Instead, we will primarily focus on cases that present particular challenges to one or more of the solution methods.

\subsection{Monomodal case: linear evaporation without coalescence}

\begin{figure}[t]
\noindent
\begin{center}
\epsfxsize=2.5in
\epsffile{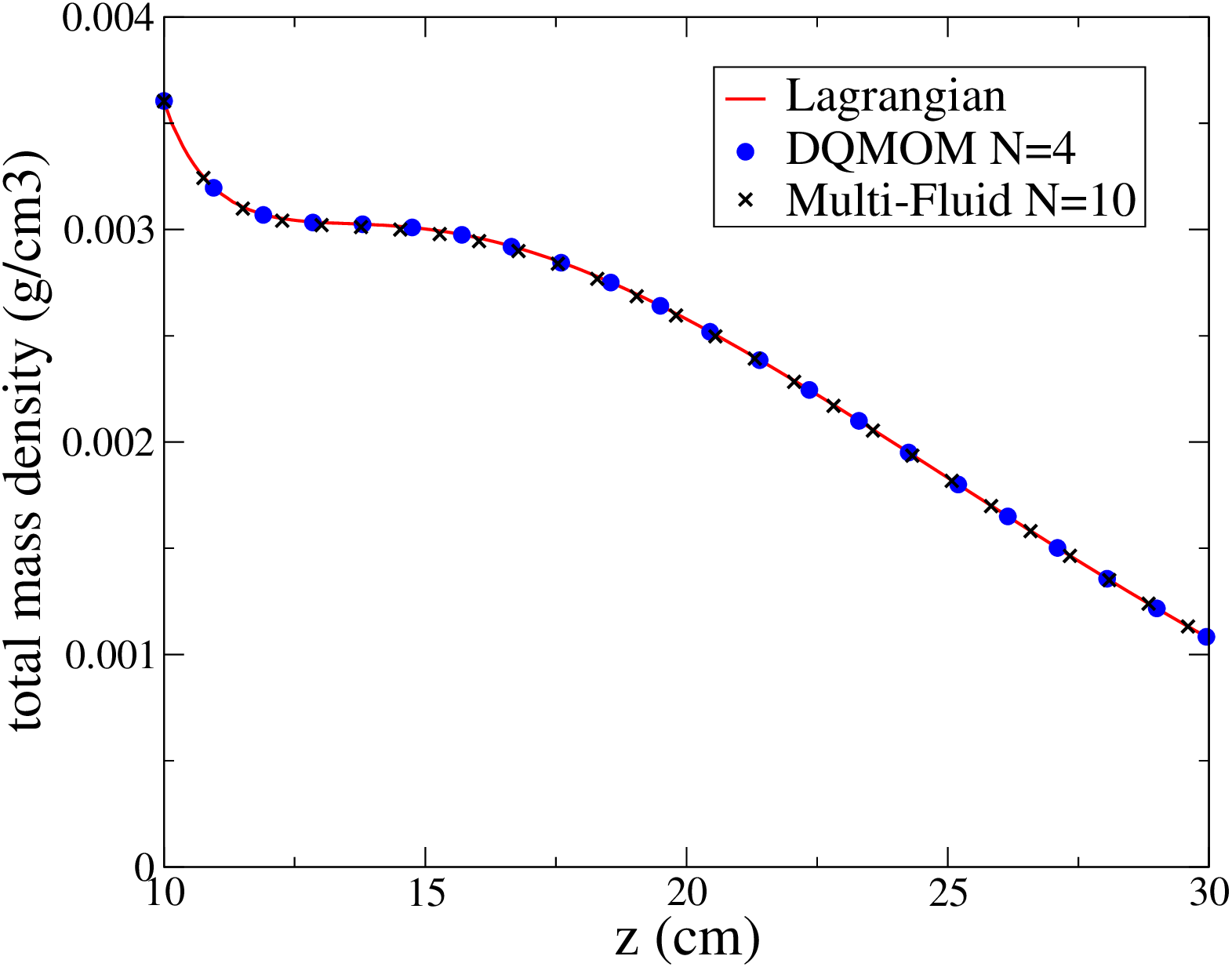} \hfill
\epsfxsize=2.5in
\epsffile{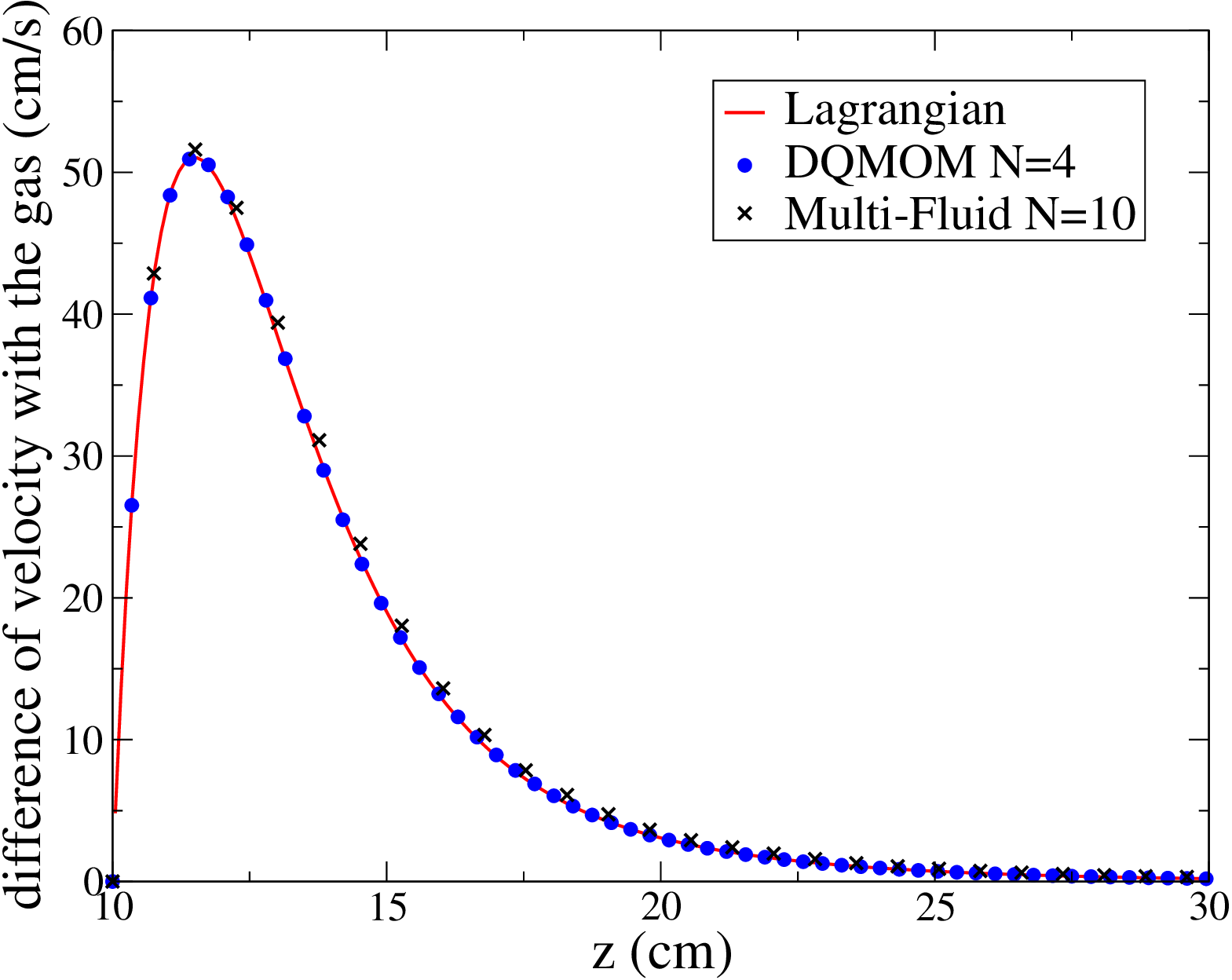}

\medskip
\epsfxsize=2.5in
\epsffile{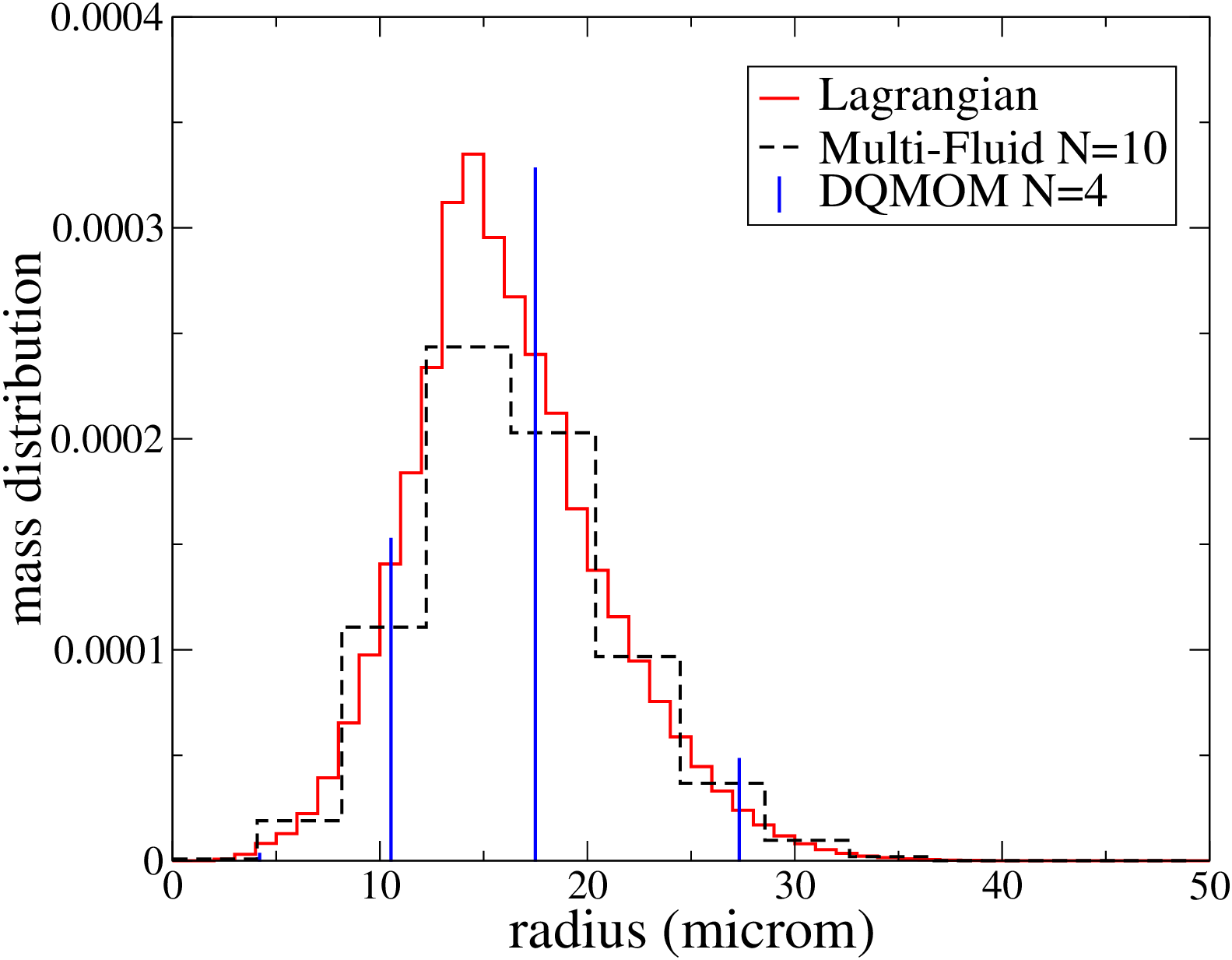} \hfill
\epsfxsize=2.5in
\epsffile{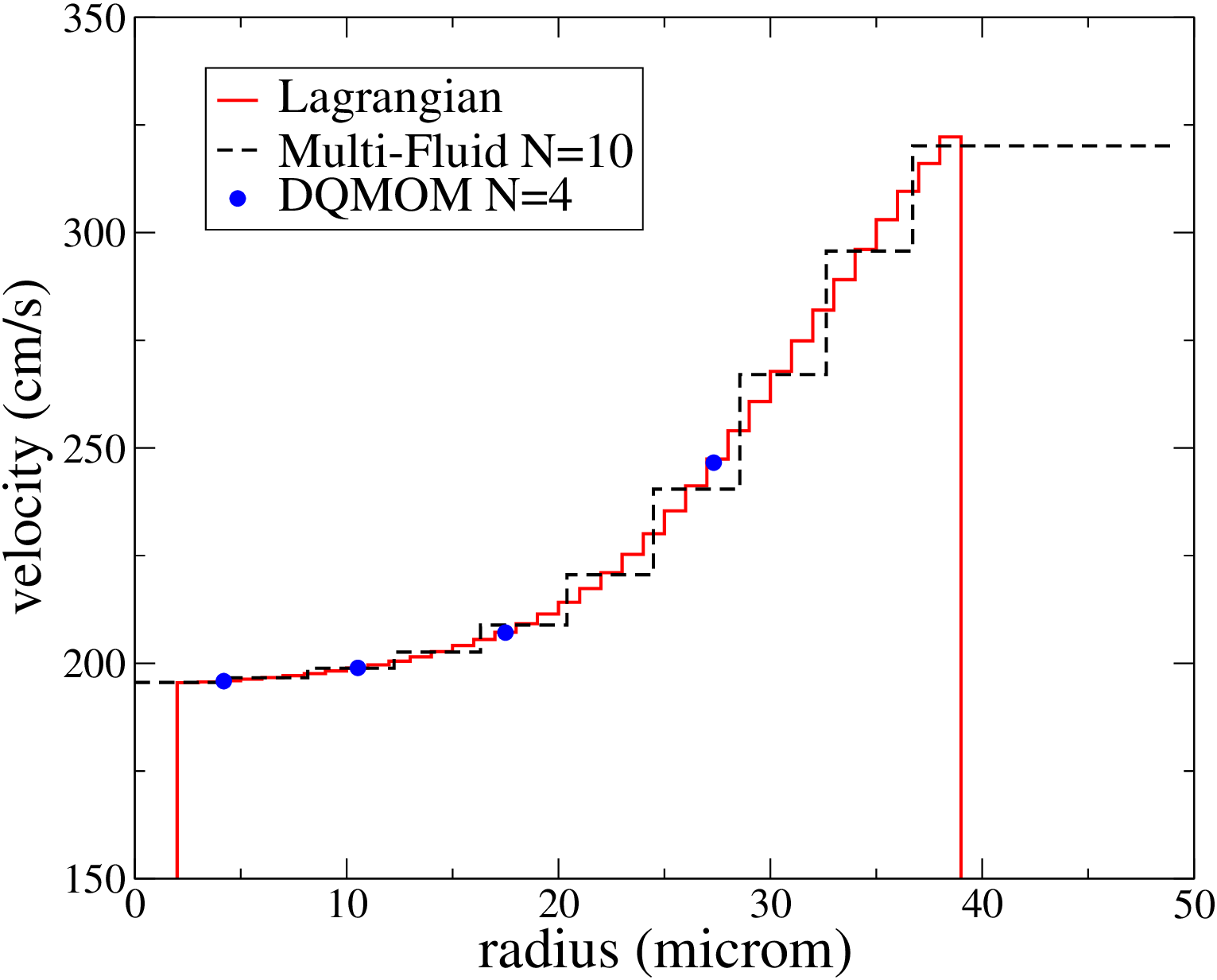}

\end{center}
\caption{Monomodal case with linear evaporation ($R_v \approx 7.126 v$). Top left: mass density. Top right: velocity difference. Bottom left: mass distribution function ($z=16$~cm). Bottom right: conditional velocity ($z=16$~cm).}
\label{fig2}
\end{figure}

We begin with a representative case where all three solution methods yield essentially identical results for all statistics.  As noted in the discussion of the methods, for linear evaporation without coalescence the DQMOM equation for each node is the same as the Lagrangian model.  Thus, the only difference between the two solution methods is that the Lagrangian method uses many more particles to represent the spray than the DQMOM method ($N=4$).  For the monomodal distribution, the multi-fluid model does not require many sections ($N=10$) to accurately capture cases with linear evaporation without coalescence.  The simulation results for the three methods are shown in Fig.~\ref{fig2}.  It can be observed that the mass density and velocity difference predicted by the three methods are nearly identical.  From the plot of the mass distribution function at $z=16$~cm, we can see that the multi-fluid model with ten sections does a good job of capturing the Lagrangian mass distribution function.  Likewise, the DQMOM weights and abscissas follow the general shape of the Lagrangian mass distribution function.  Finally, for the conditional velocity $\langle u_z | r \rangle$ we see that all three methods produce the same curve. 
We should note that for cases without coalescence the Lagrangian simulations predict essentially no velocity dispersion about the conditional value.  In other words, conditional velocity fluctuations defined by $u'(r) \equiv \langle (u_z - \langle u_z | r \rangle)^2 | r \rangle^{1/2}$ are null.  This is exactly one of the necessary conditions evoked when deriving the multi-fluid model, which would explain why its predictions for this case are in excellent agreement with the Lagrangian method.

\subsection{Monomodal case: nonlinear evaporation without coalescence}

Cases with nonlinear evaporation result in a loss of droplets in finite time, which translates into a nonzero flux $\psi (t)$ in DQMOM. For the monomodal case without coalescence, we expect the flux term to be a smooth function of $t$, and thus it cannot be neglected.  In the multi-fluid model, the flux is computed directly from the shape of the first section (i.e., the section near the origin)
and does not yield any difficulty.  

\begin{figure}[t]
\noindent
\begin{center}

\epsfxsize=2.5in
\epsffile{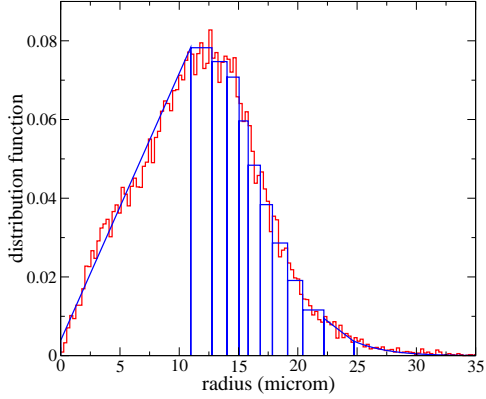}

\end{center}
\caption{Monomodal distribution function with optimal sections.}
\label{fig3}
\end{figure}

\begin{figure}[t]
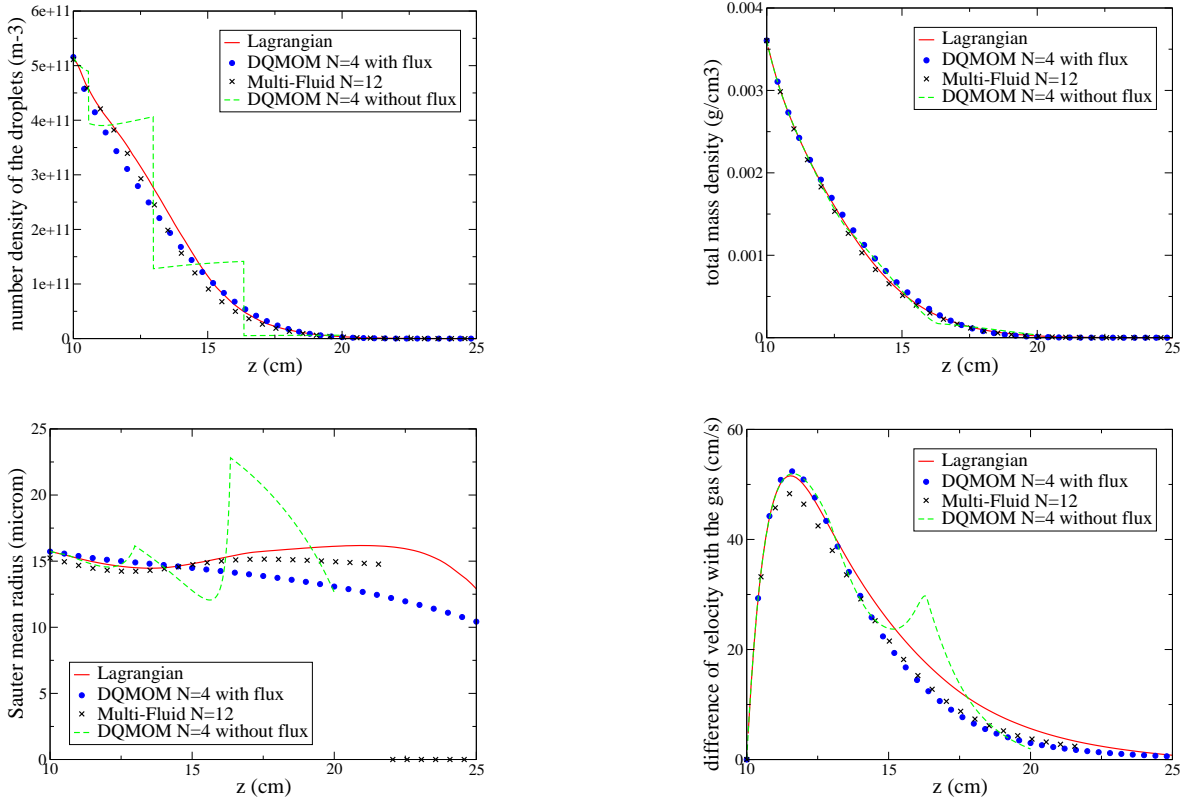

\noindent
\begin{center}

\epsfxsize=2.5in
\epsffile{./Figures_Article/comp_number_Rscst_nonrafine.eps} \hfill
\epsfxsize=2.5in
\epsffile{./Figures_Article/comp_mass_Rscst_nonrafine.eps}

\bigskip
\epsfxsize=2.5in
\epsffile{./Figures_Article/comp_rSauter_Rscst_nonrafine.eps} \hfill
\epsfxsize=2.5in
\epsffile{./Figures_Article/comp_vita_Rscst_nonrafine.eps}

\end{center}
\caption{Monomodal case with nonlinear evaporation. Top left: number density. Top right: mass density. Bottom left: Sauter mean radius. Bottom right: velocity difference.}
\label{fig4}
\end{figure}

In our multi-fluid simulations, we use the ``optimal'' choice of sections with $N=12$ shown in Fig.~\ref{fig3} \cite{lsgsm04}. Obviously, a finer discretization (larger $N$) could be used in the multi-fluid model to attain closer agreement with the Lagrangian method, but this would increase the computational cost. Note that the first section is represented by a constant slope, which corresponds to a constant flux level at each time step.  For the DQMOM, we use $N=4$ and the evaporative flux is computed using the ratio constraints introduced in Section~\ref{ef}. It can be noticed that the increase of $N$ do not imply an increase of the number of conserved moments during the evaporation step since the number of ratio constraints is also increasing in the same way. The value of $N$ is then conditioned by the capacity of the method to follow the dynamics of droplets of different sizes.
Representative results for the three solution methods are shown in Fig.~\ref{fig4}. In general, all three methods produce very similar predictions.
From the number density, we can observe that DQMOM with the ratio constraints does a good job of predicting the loss of droplets due to evaporation. 
Likewise, the mass densities found from all methods are very close.
We should note that for $z>20$~cm the number of remaining droplets is very small and the statistics computed from the Lagrangian method are subject to statistical errors.
Comparing the Sauter mean radii predicted by the three methods, we can observe that the agreement is generally satisfactory up to $z=20$~cm.
The DQMOM shows the largest deviation from the Lagrangian method at $z=20$~cm due to errors in the flux model, but the agreement is still acceptable.
The differences in the Sauter mean radius are reflected in the predictions of the velocity difference.
In general, droplets with a larger radius will have a higher velocity difference.
Thus, we see that initially the Sauter mean radius predicted by DQMOM is larger than that from the Lagrangian method, resulting in a slightly higher velocity difference at $z=12$~cm.  Later on ($z>15$~cm) this trend is reversed. Finally, we can note that neglecting the flux in DQMOM yields poor predictions of number density 
since we can observe the artificial jumps in the number density related to 
the singular fluxes associated to one abscissa crossing the zero size limit,
as well as the oscillating dynamics of 
the Sauter mean radius for this case.

\subsection{Bimodal case: nonlinear evaporation without coalescence}

\begin{figure}[t]
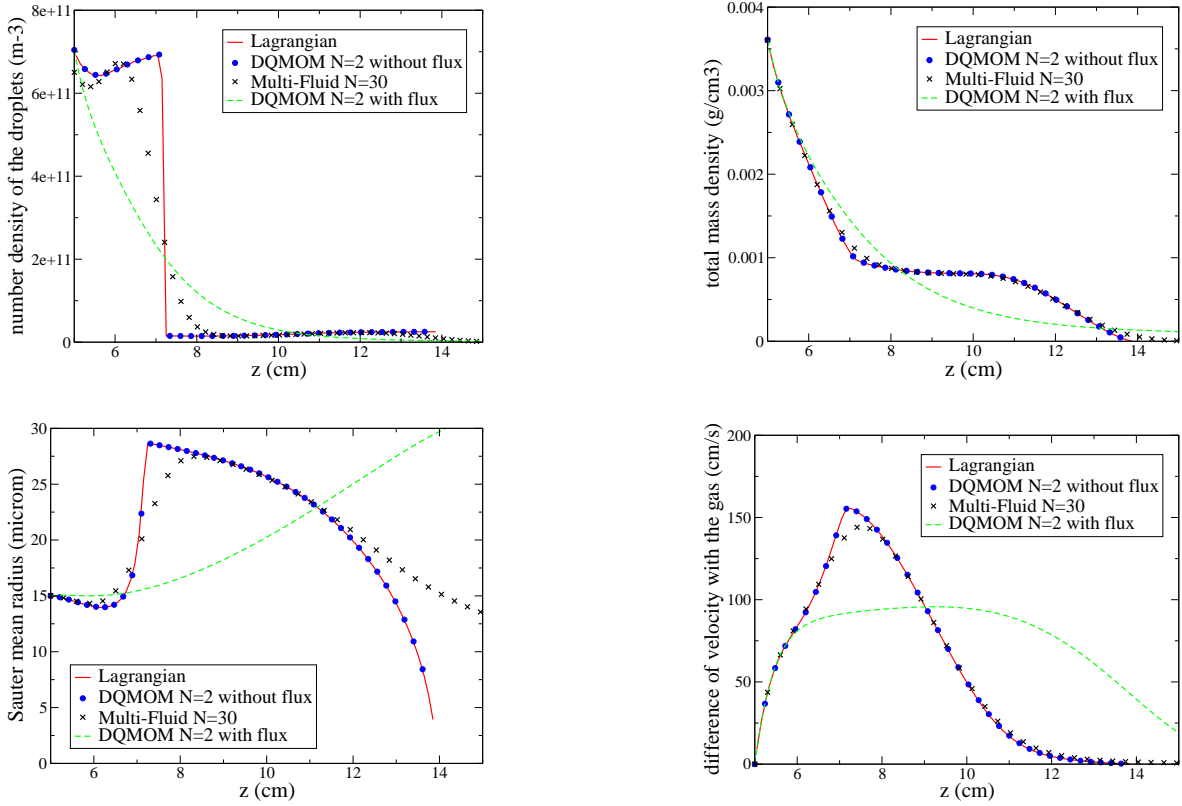

\noindent
\begin{center}

\epsfxsize=2.5in
\epsffile{./Figures_Article/comp_number_bi_Rscst_nonrafine.eps} \hfill
\epsfxsize=2.5in
\epsffile{./Figures_Article/comp_mass_bi_Rscst_nonrafine.eps}

\bigskip
\epsfxsize=2.5in
\epsffile{./Figures_Article/comp_rSauter_bi_Rscst_nonrafine.eps} \hfill
\epsfxsize=2.5in
\epsffile{./Figures_Article/comp_vita_bi_Rscst_nonrafine.eps}

\end{center}
\caption{Bimodal case with nonlinear evaporation. Top left: number density. Top right: mass density. Bottom left: Sauter mean radius. Bottom right: velocity difference.}
\label{fig5}
\end{figure}

By changing from the monomodal to the bimodal distribution, we change the nature of the initial distribution function and thus the nature of the numerical difficulties.
For the multi-fluid model, the bimodal case is difficult because a relatively large number of sections ($N=30$) is needed to capture the two peaks with acceptable numerical diffusion.
The use of a second-order method developed in \cite{laurent06} would reduce this number to around 10; however, it would still be difficult to describe Dirac delta function by a finite-volume approximation.
On the other hand, this case is ``optimal'' for DQMOM because only two ($N=2$) abscissas are required (one for each peak) and the flux is null, expect when a peak passes the origin.
In Fig.~\ref{fig5} results from the three simulation methods are shown and it is clear that DQMOM performs extremely well for this case by setting $\psi=0$.
For example, the number density function shows step changes at $z=7.2$~cm and $13.8$~cm (i.e., when a peak passes the origin), and DQMOM exactly reproduces this behavior.
With $N=30$, the multi-fluid model does a good job of predicting the mass density.
However, from the plots of number density and Sauter mean radius, we can observe the negative effects of numerical diffusion, which tends to smooth out the peaks in the distribution (the method has been shown to be first order in the droplet size discretization step in \cite{laurent06}, where some higher-order methods have been proposed).
Nevertheless, all three methods yield reasonable predictions for all of the cases without coalescence.
We should note, however, that for more complicated initial distributions (e.g., delta functions combined with smooth functions) specifying the evaporative flux in DQMOM may be problematic.

\subsection{Monomodal case: linear evaporation with coalescence}

\begin{figure}[t]
\noindent
\begin{center}

\epsfxsize=2.5in
\epsffile{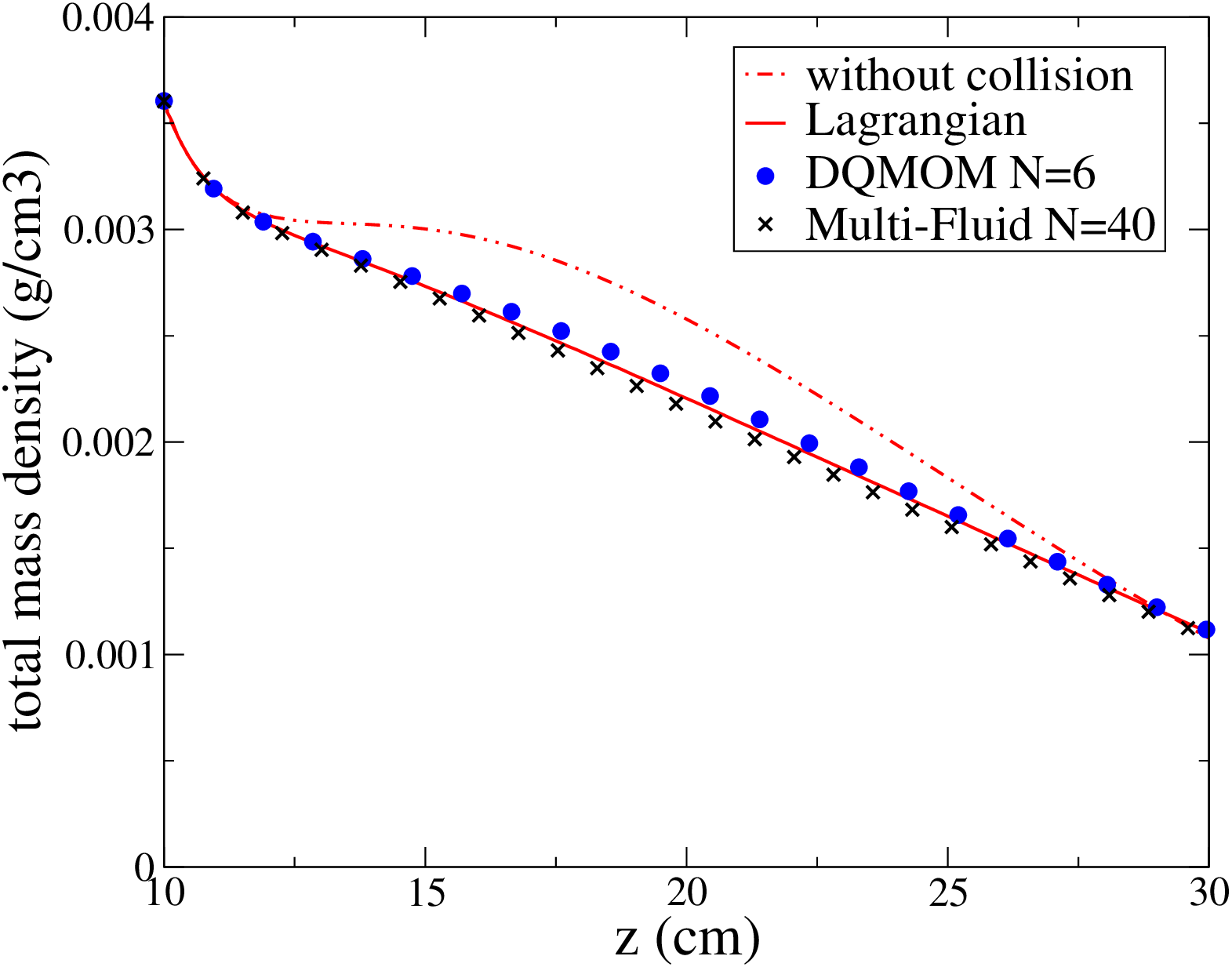} \hfill
\epsfxsize=2.5in
\epsffile{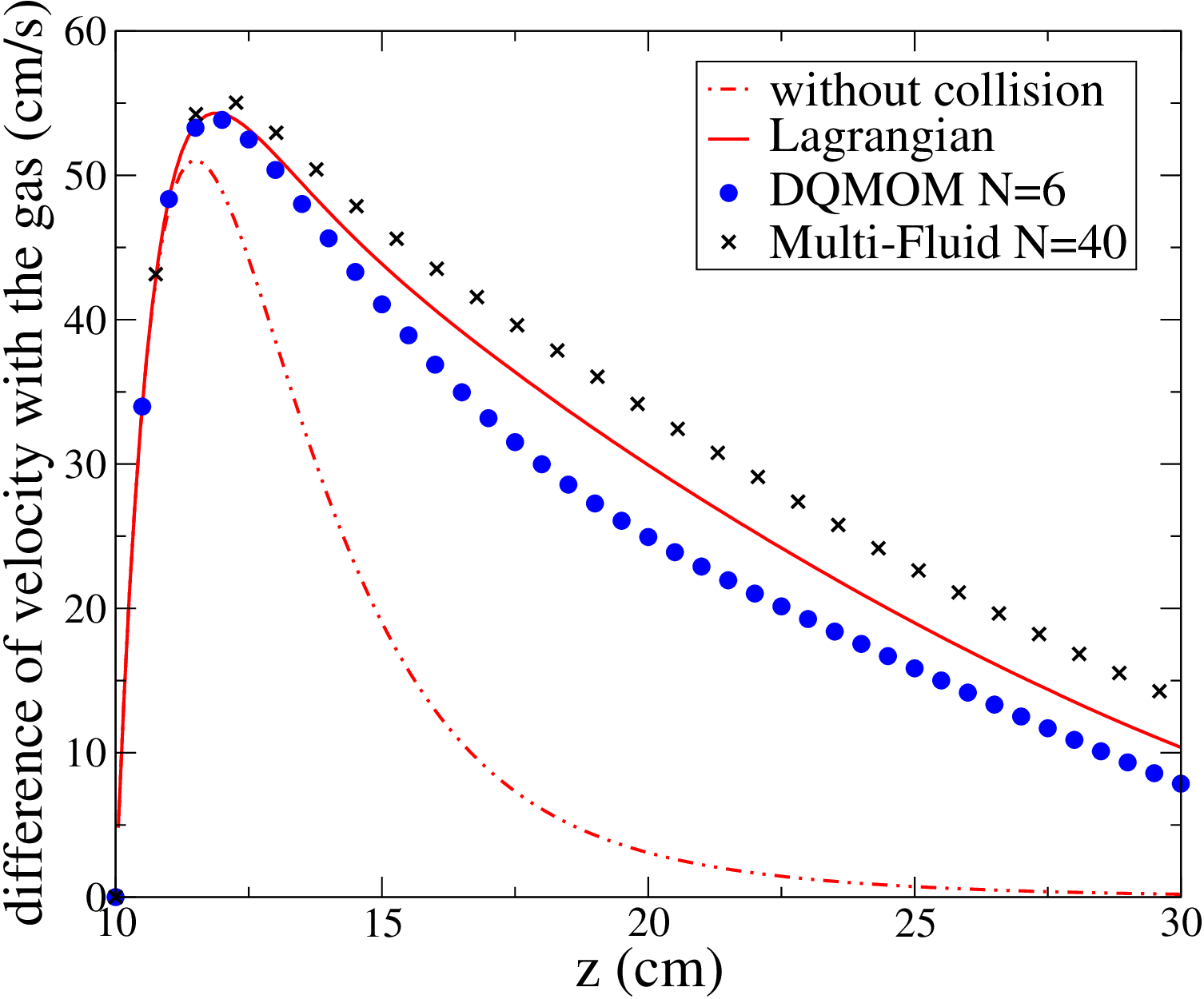}

\bigskip
\epsfxsize=2.5in
\epsffile{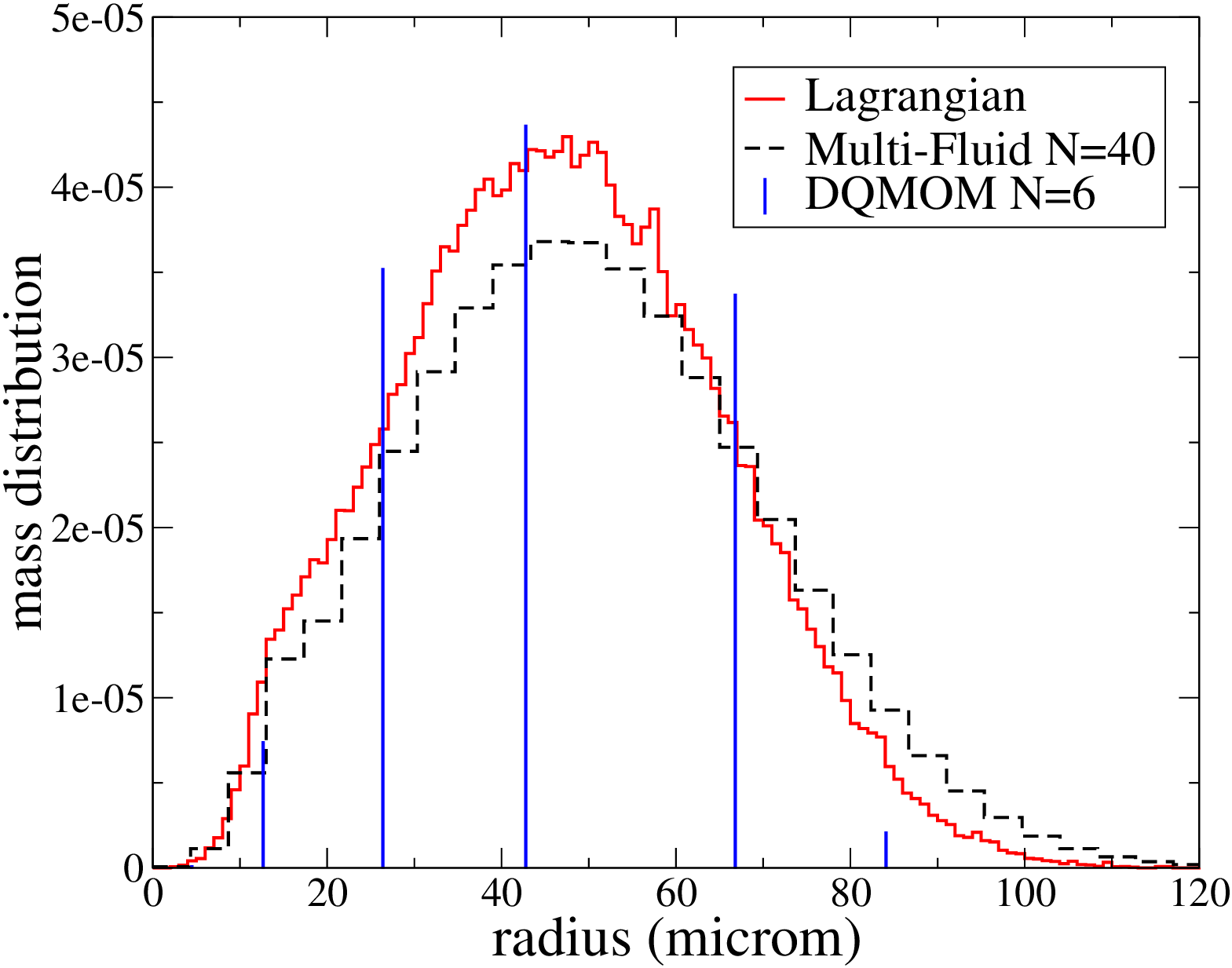} \hfill
\epsfxsize=2.5in
\epsffile{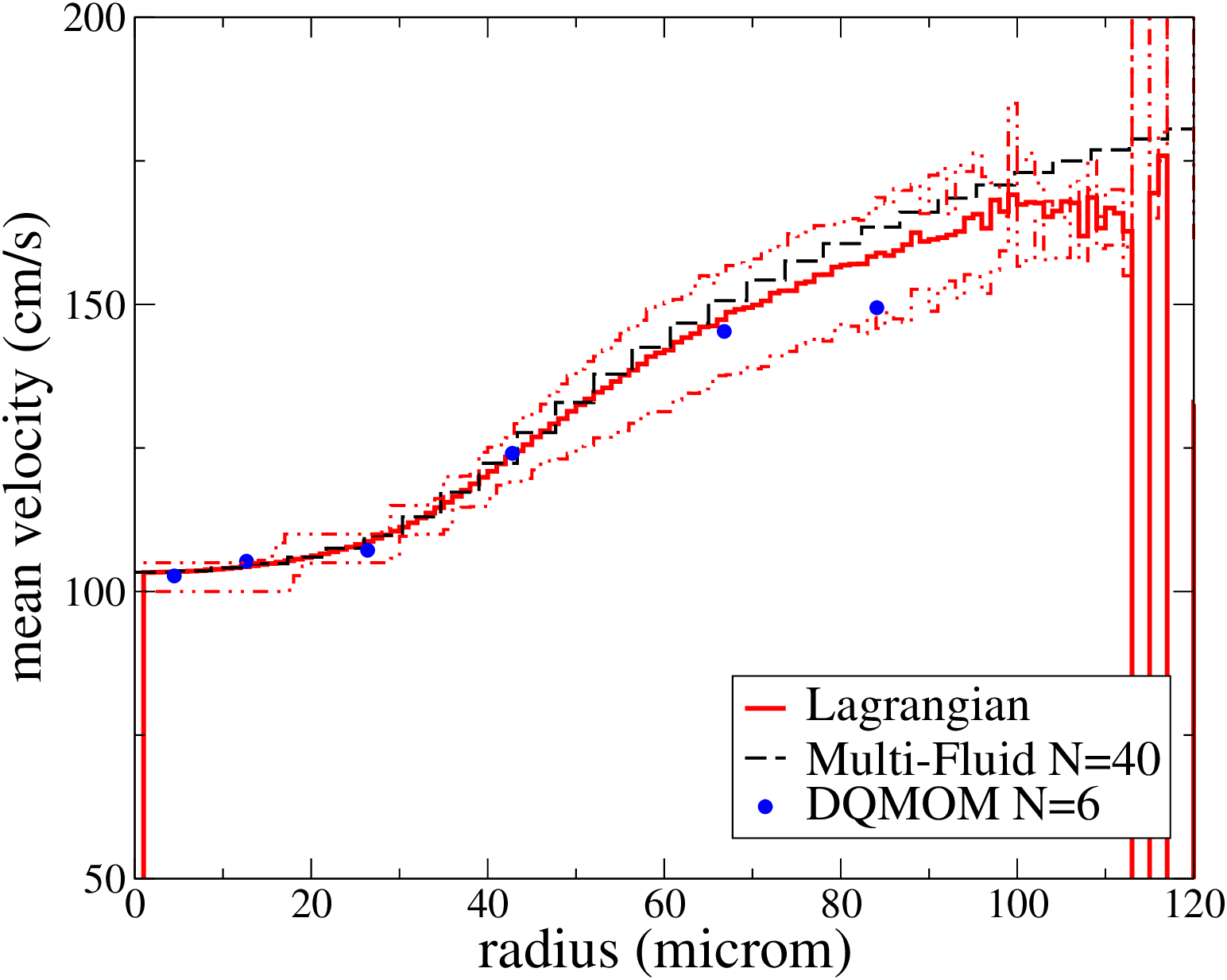}

\end{center}
\caption{Monomodal case with linear evaporation ($R_v \approx 7.126 v$) and coalescence. Top left: mass density. Top right: velocity difference. Bottom left: mass distribution function ($z=22$~cm). Bottom right: conditional velocity ($z=22$~cm).}
\label{fig6}
\end{figure}

We now turn to the more difficult cases that include coalescence.
As mentioned earlier, the coalescence of droplets with different volumes (and velocities) will lead to velocity dispersion ($u'(r)>0$).
Physically, this implies that two droplets with the same volume will have a nonzero probability of colliding (due to the difference in velocity).
Thus, the rate of coalescence when $u'(r)>0$ will be larger than when $u'(r)=0$.
Numerical approximations (such as the multi-fluid model) that assume $u'(r)=0$ should therefore predict smaller droplets than the Lagrangian method.  In Fig.~\ref{fig6} we present results for the three methods for linear evaporation ($\psi=0$) with coalescence.
>From the velocity difference we can observe that coalescence leads to a slower relaxation to the gas velocity due to formation of larger droplets than without coalescence.
Note that in general all three methods predict similar results for the velocity difference.
However, due the differences in the predictions of the shape of the mass distribution function, the multi-fluid model predicts slightly slower relaxation and the DQMOM slightly faster than is found with the Lagrangian method.
Comparing with Fig.~\ref{fig1}, we can observe that coalescence leads to much larger droplets than are present in the initial distribution function. In general, the multi-fluid model predicts a slightly larger number of droplets above $80$~$\mu$m than the Lagrangian method.
Nevertheless, the predictions are in reasonably good agreement.
The predictions for the conditional velocity $\langle u_z | r \rangle$ are also good.
Finally, note that we used $N=6$ with DQMOM, the reason for which will be discussed for a more difficult case in Section~\ref{ssevap}.

\subsection{Bimodal case: linear evaporation with coalescence}

\begin{figure}[t]
\noindent
\begin{center}

\epsfxsize=2.5in
\epsffile{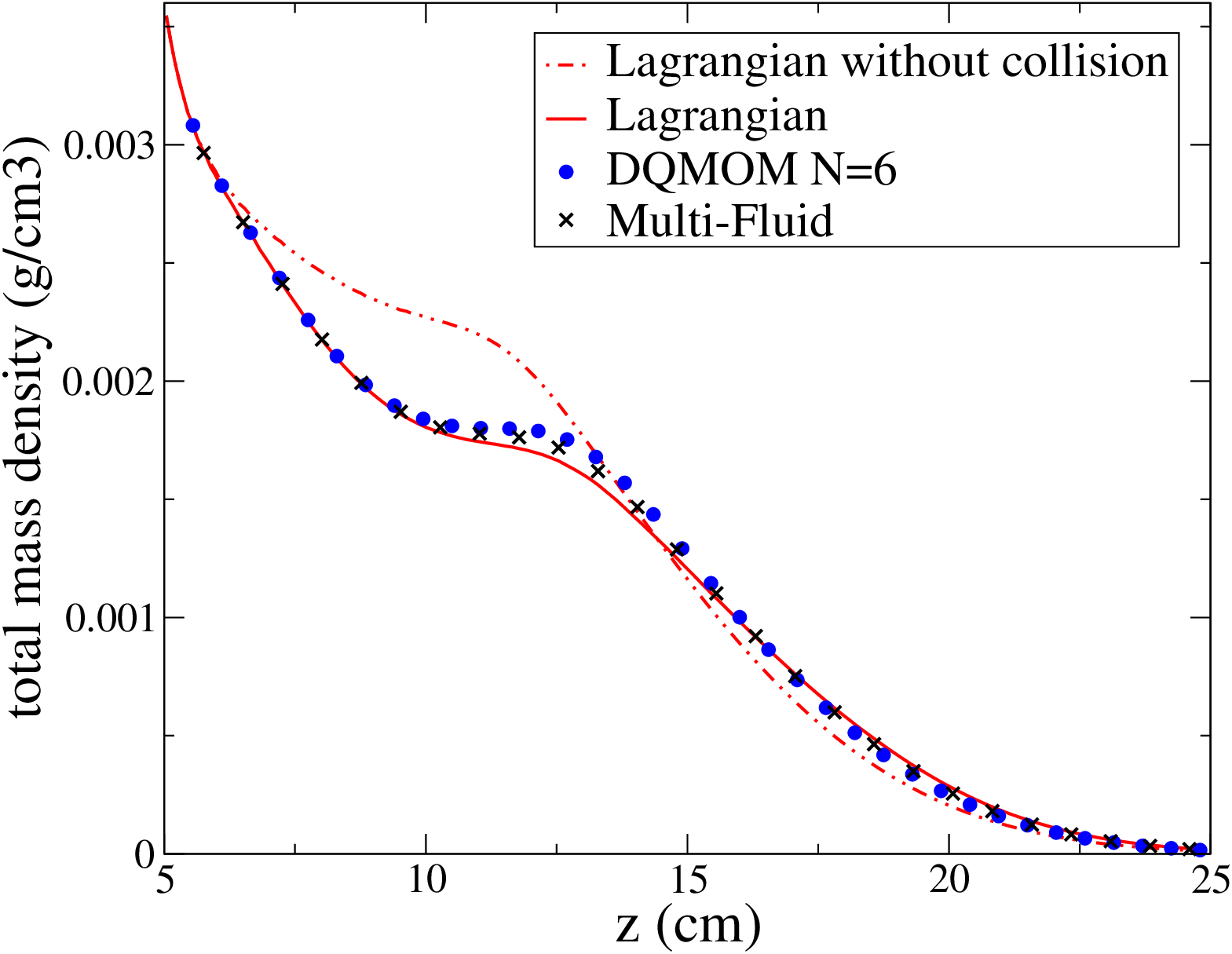} \hfill
\epsfxsize=2.5in
\epsffile{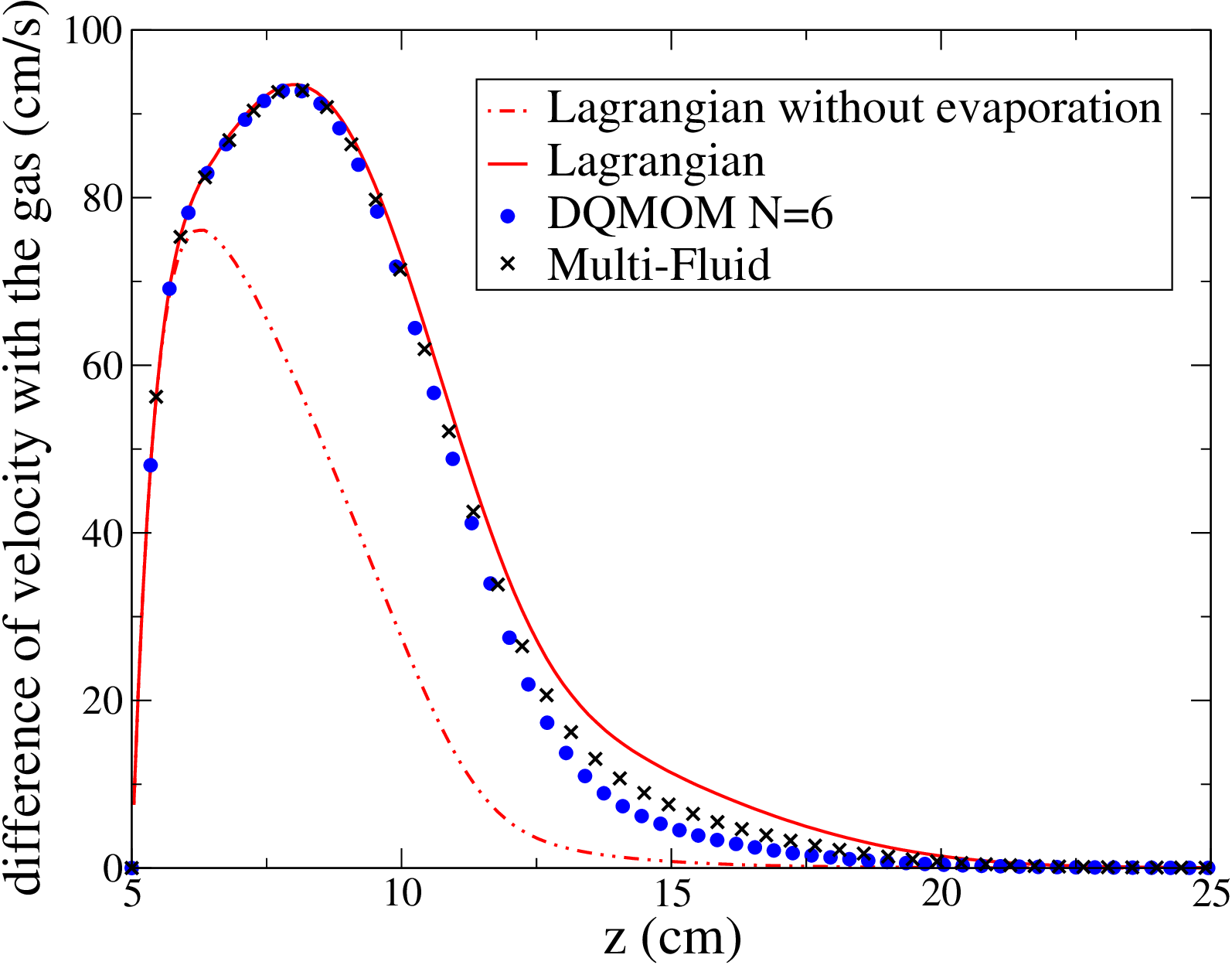}

\bigskip
\epsfxsize=2.5in
\epsffile{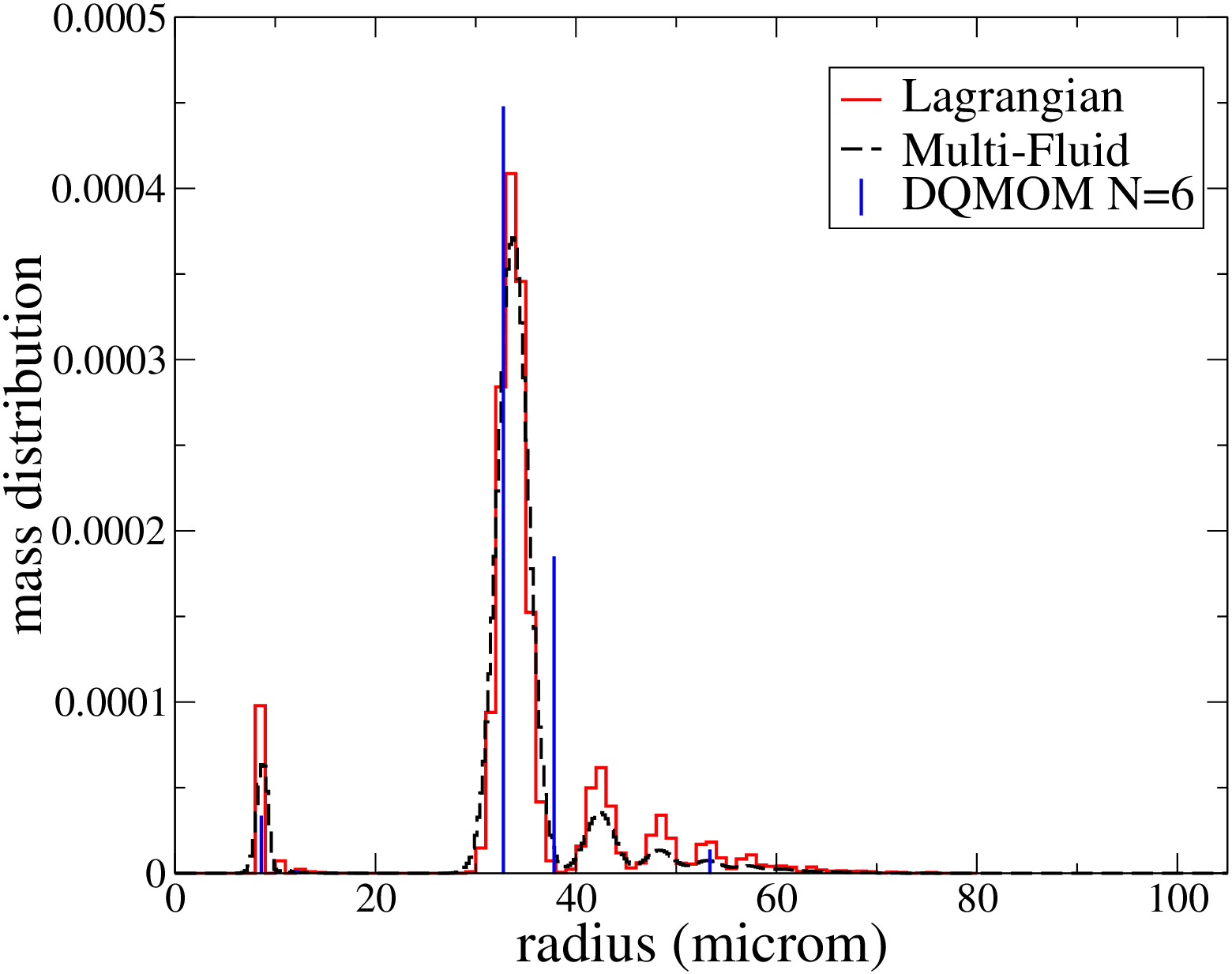} \hfill
\epsfxsize=2.5in
\epsffile{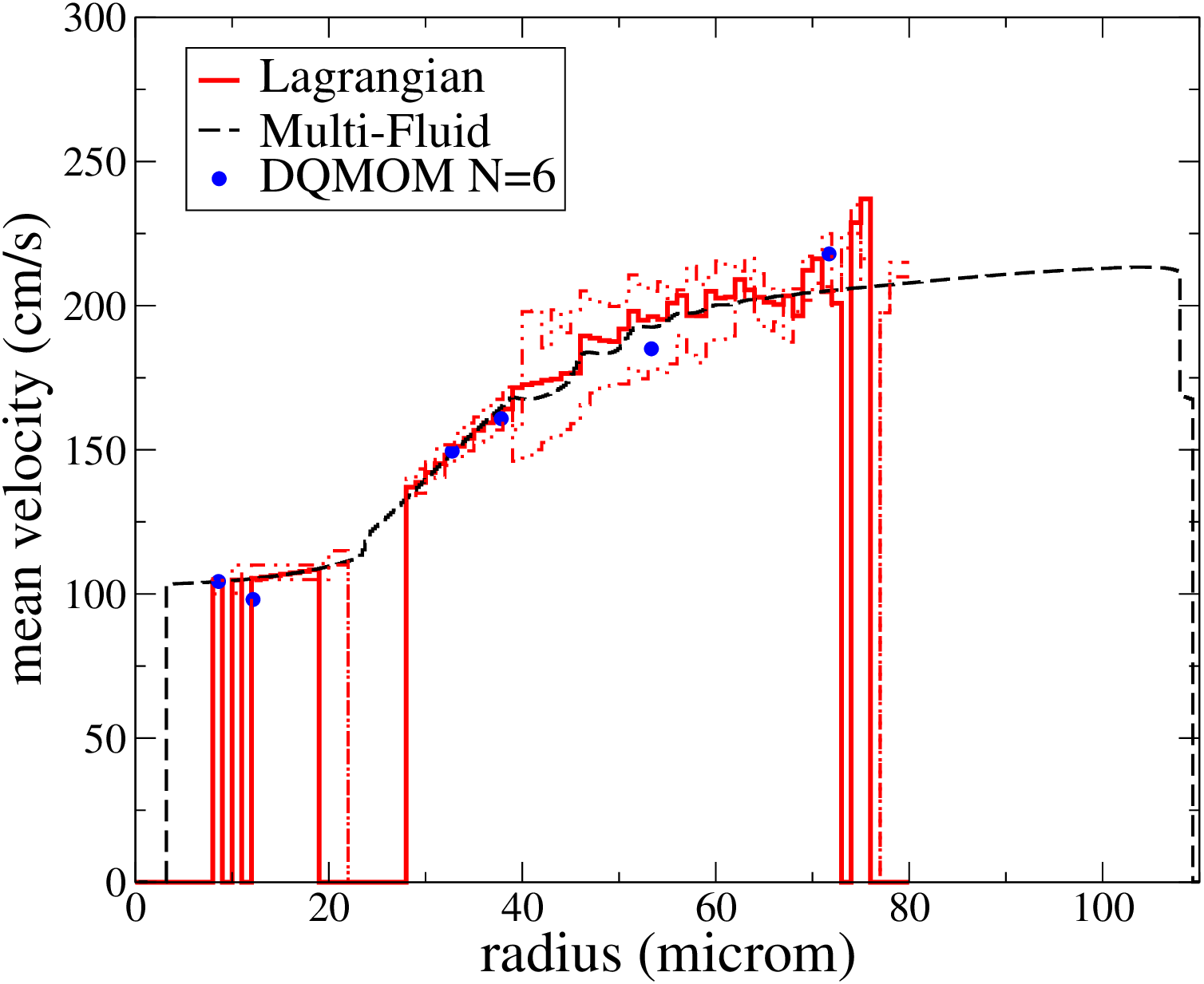}

\end{center}
\caption{Bimodal case with linear evaporation ($R_v \approx 14.252 v$) and coalescence. Top left: mass density. Top right: velocity difference. Bottom left: mass distribution function ($z=11$~cm). Bottom right: conditional velocity ($z=11$~cm).}
\label{fig7}
\end{figure}

We now consider a more difficult case where the initial distribution is bimodal.  As discussed previously, the peaks in the distribution are difficult to resolve accurately in the multi-fluid model
with a limited number of
sections
. When combined with coalescence, this has important consequences because numerical
diffusion
can lead to spurious coalescence of droplets with slightly different volumes (and hence velocities)
as observed in \cite{lmv04}.
For example, with the bimodal distribution with droplets of radii $10$ and $30$~$\mu$m, coalescence cannot produce droplets below $30$~$\mu$m.  However, spurious coalescence between droplets of radii near $10$~$\mu$m leads to droplets in the range below $30$~$\mu$m. We overcome this difficultly by using a large number of sections ($N=500$) in the multi-fluid model.
This number could also be reduced
by using a second-order method for the evaporation such as the one of \cite{laurent06} but this is not the point we want to make with this configuration.
Note that the same problem arises in the Lagrangian method when the spatial cell size $\Delta z$ is too large.  While DQMOM does not suffer from spurious coalescence, the bimodal case is still difficult because the initially two-peak distribution will quickly form multiple peaks due to pair-wise collisions. With $N=6$ in DQMOM, it is at best possible to represent six peaks.  Results for the three methods are shown in Fig.~\ref{fig7} where it can be seen that the mass density and the velocity difference are reasonably well predicted by the multi-fluid model and DQMOM.  From the mass distribution function at $z=11$~cm, the multi-peak structure due to coalescence is quite evident, as is the slight numerical diffusion 
in the multi-fluid model (even with $N=500$, but this is expected since this is a first-order method).  Note that DQMOM with $N=6$ has two abscissas at points corresponding to the major peaks ($10$ and $33$~$\mu$m), and the remaining abscissas at points without major peaks.  Comparisons of the conditional velocity predicted by the three methods are also quite favorable for this difficult case.

\subsection{Monomodal case: nonlinear evaporation with coalescence}

\begin{figure}[t]
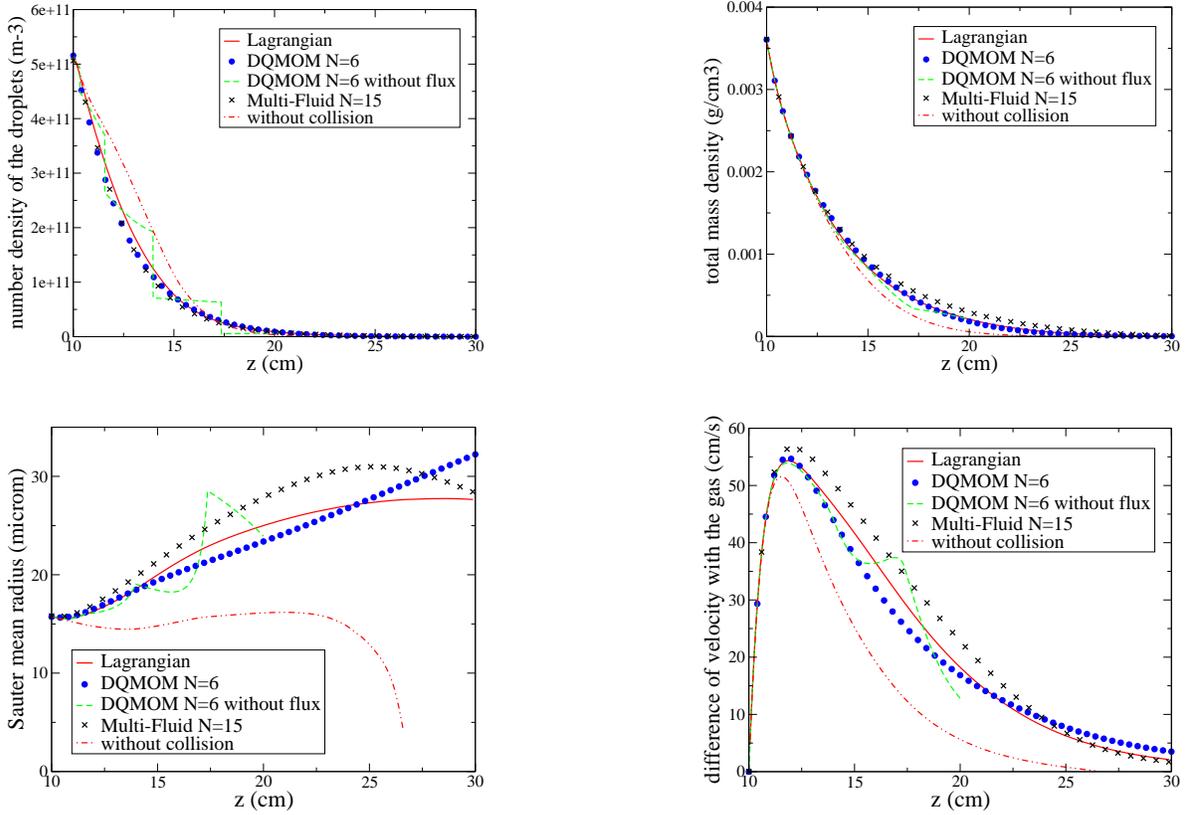

\noindent
\begin{center}

\epsfxsize=2.5in
\epsffile{./Figures_Article/comp_number_Rscst_col_nonrafine.eps} \hfill
\epsfxsize=2.5in
\epsffile{./Figures_Article/comp_mass_Rscst_col_nonrafine.eps}

\bigskip
\epsfxsize=2.5in
\epsffile{./Figures_Article/comp_rSauter_Rscst_col_nonrafine.eps} \hfill
\epsfxsize=2.5in
\epsffile{./Figures_Article/comp_vita_Rscst_col_nonrafine.eps}

\end{center}
\caption{Monomodal case with nonlinear evaporation and coalescence. Top left: number density. Top right: mass density. Bottom left: Sauter mean radius. Bottom right: velocity difference.}
\label{fig8}
\end{figure}

We will now consider the more physically relevant case of nonlinear evaporation.
As discussed earlier, the evaporative flux for this case is nonzero, so we will need to model it in DQMOM.  Here, we consider two models for $\psi$: (a) ratio constraints and (b) $\psi=0$.
Because the initial distribution is monomodal, we might expect that using ratio constraints is always a better choice.
On the other hand, if coalescence is much faster than evaporation, it might happen that droplets grow faster than they disappear so that the evaporative flux is closer to zero.
For the multi-fluid model, we use $N=15$ sections. Results for the three methods are shown in Fig.~\ref{fig8}.
The number density illustrates the effect of the choice of $\psi$ in DQMOM.
With $\psi=0$, the number density changes discontinuously whenever an abscissa passes the origin.
However, DQMOM with ratio constraints yields predictions very similar to the other two methods.
Likewise, the mass density is predicted to be very similar for all three methods; however, using zero flux with DQMOM is slightly worse.
The predictions for the Sauter mean radius show opposing trends.
In general, the multi-fluid model overpredicts the mean radius (i.e. predicts too much coalescence), while DQMOM underpredicts it.
As before, for the DQMOM predictions, the results with the ratio constraints are best.
The predictions for the velocity difference follows the same trend.
As discussed in the next example, the differences observed between the Lagrangian method and the two Eulerian methods is likely due to the latter's inability to capture velocity dispersion.
Moreover, we have used $N=6$ with DQMOM since, as shown in Section~\ref{ssevap}, it is adequate to describe coalescence phenomenon for this particular set of moments.

\subsection{Monomodal case: coalescence with no evaporation}
\label{ssevap}

\begin{figure}[t]
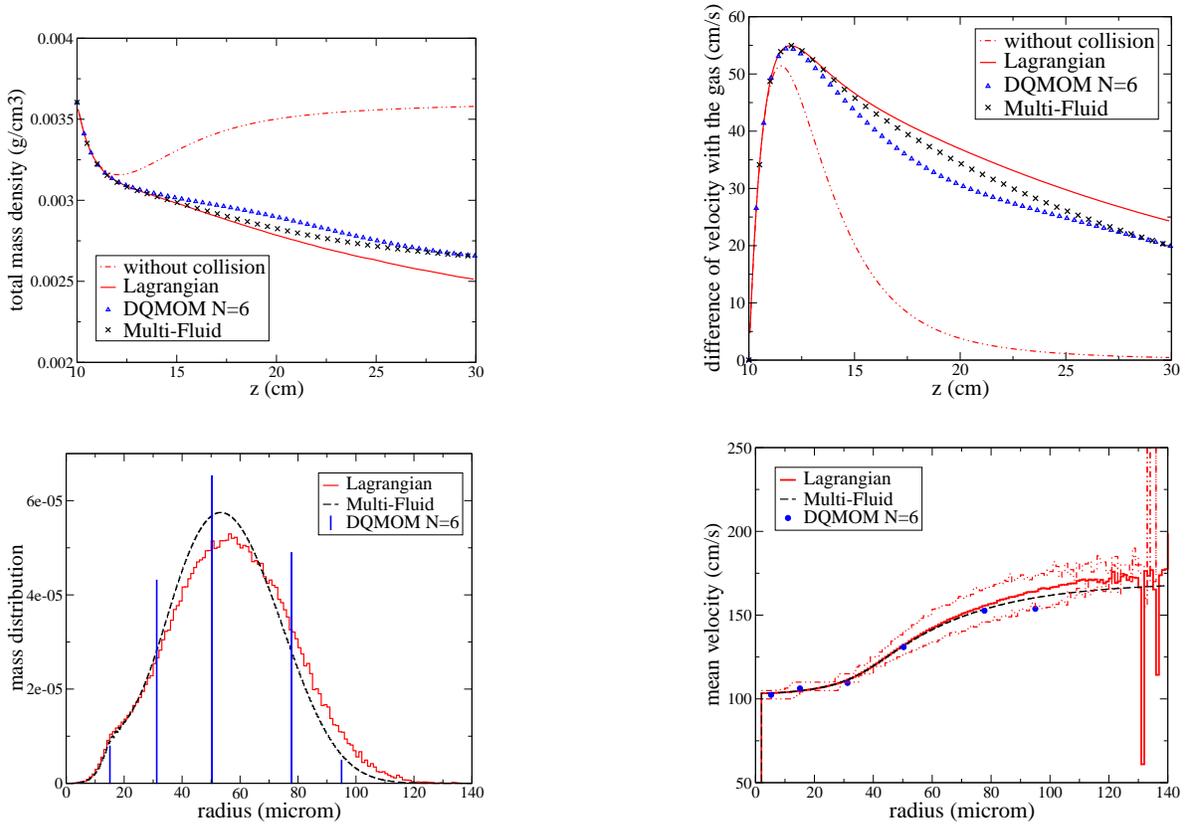

\noindent
\begin{center}

\epsfxsize=2.5in
\epsffile{./Figures_Article/comp_mass_ssevap_rafine.eps} \hfill
\epsfxsize=2.5in
\epsffile{./Figures_Article/comp_vita_ssevap_rafine.eps}

\bigskip
\epsfxsize=2.5in
\epsffile{./Figures_Article/distrm_ssevap_col_z22.eps} \hfill
\epsfxsize=2.5in
\epsffile{./Figures_Article/distrv_ssevap_col_z22.eps}

\end{center}
\caption{Monomodal case with coalescence but no evaporation. Top left: mass density. Top right: velocity difference.
Bottom left: mass distribution ($z=22$~cm). Bottom right: conditional velocity ($z=22$~cm).}
\label{fig9}
\end{figure}

\begin{figure}[t]
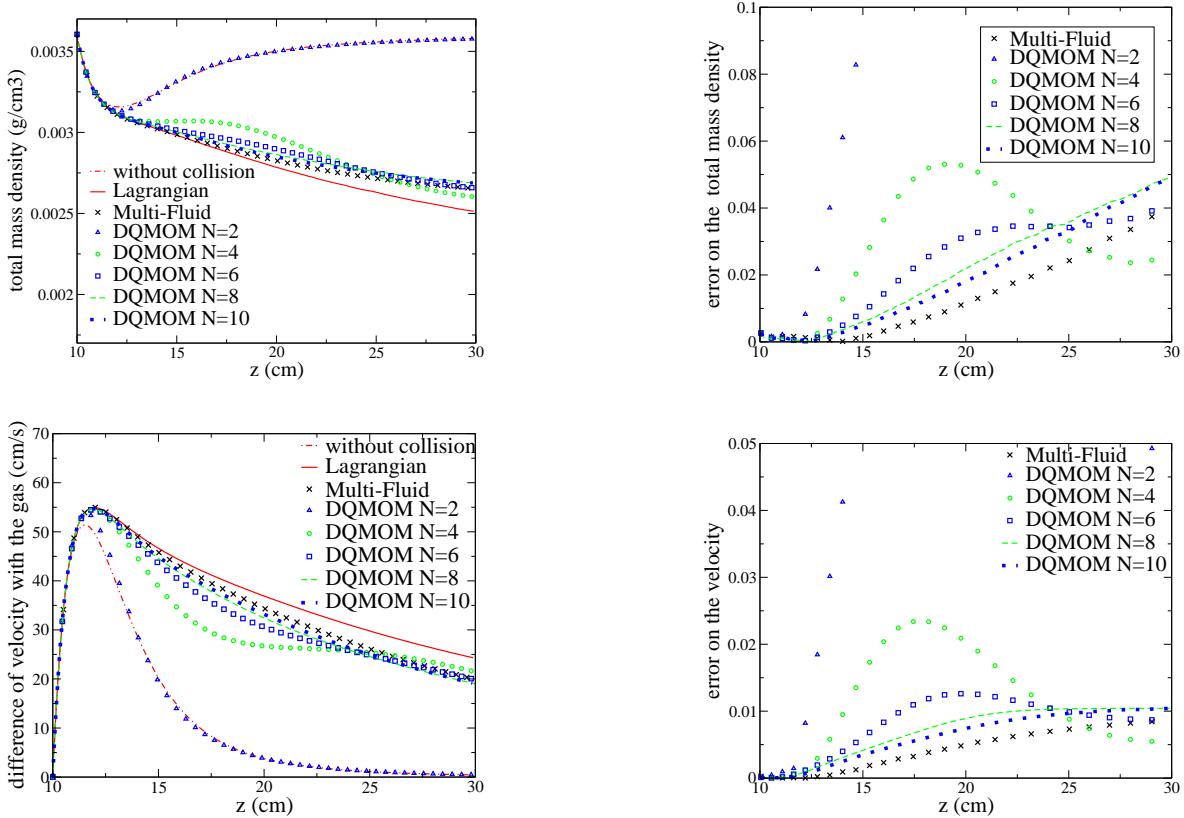

\noindent
\begin{center}

\epsfxsize=2.5in
\epsffile{./Figures_Article/comp_mass_ssevap.eps} \hfill
\epsfxsize=2.5in
\epsffile{./Figures_Article/comp_dmass_ssevap.eps}

\bigskip
\epsfxsize=2.5in
\epsffile{./Figures_Article/comp_vita_ssevap.eps}\hfill
\epsfxsize=2.5in
\epsffile{./Figures_Article/comp_dvita_ssevap.eps}

\end{center}
\caption{Monomodal case with coalescence but no evaporation.
Top left: mass density. Top right: error on the mass density.
Bottom left: velocity difference. Bottom right: error on the velocity difference.}
\label{fig89}
\end{figure}

In order to highlight the role of coalescence on determining the evolution of the number density function, we now consider a case with no evaporation.  For this case, droplets will grow continuously due to coalescence, and velocity dispersion will enhance the collision rate and lead to even larger droplets.  Because the multi-fluid model uses fixed sections, it is necessary to fix the maximum radius at $200$~$\mu$m with $N=500$ in order to capture the largest droplets at $z=30$~cm.  In contrast, the abscissas in DQMOM move in phase space to accommodate growth.  Nevertheless, we can anticipate that the number of abscissas will affect the accuracy of the DQMOM predictions.
In Fig.~\ref{fig89} it can be observed that when the number of moments increases, the accuracy of the DQMOM solution increases, from something almost ignoring the coalescence phenomenon with $N=2$ to a saturation of the accuracy for $N\ge 8$. 
Indeed, the accuracy of the DQMOM for the description of the coalescence is related to the accuracy of the approximation of the coalescence operator by the quadrature formula (\ref{cmu3}) and (\ref{cpu3}), which increases with $N$.
Since the results are quite good and at a low cost (and the linear system is reasonably well conditioned), we will use $N=6$ for comparisons with the other two methods.

As mentioned earlier, without coalescence or evaporation all three methods predict essentially identical results.  In Fig.~\ref{fig9} the results for the pure coalescence case are shown.  Notice that the mass density does not decrease to zero because there is no evaporation; however, it does change due to transport.  From the velocity difference, we can see that the multi-fluid model and DQMOM overpredict the relaxation rate.  As discussed previously, this is due to both methods underpredicting the mean droplet size.  From the mass distribution functions at $z=22$~cm, we can observe that the Lagrangian method has more droplets with radii above $80$~$\mu$m than the multi-fluid model, which is consistent with the observed trend in the velocity difference. In order to explore the link between velocity dispersion and coalescence, we have computed 50\% probability intervals for the conditional velocity.  These are defined to be the values of $v$ for which the conditional velocity PDF $f(v|r)$ is fifty percent of its peak value. Note that in the absence of velocity dispersion $f(v|r)$ is a delta function centered at $\langle u_z | r \rangle$, so the width of the intervals is a measure of dispersion.
From the plot of conditional velocity, we can note that for large droplets the velocity dispersion is significant.
We can also note that using DQMOM
essentially results in
points along the curve $\langle u_z | r \rangle$, i.e., increasing $N$ with the same choice of moments does not capture the velocity dispersion. 


\section{Conclusions}

In this work, we have implemented DQMOM to treat the Williams spray equation that describes evaporation, acceleration and coalescence of liquid droplets in a laminar gas flow.
The derivation of the DQMOM equations was shown to be a straightforward task, and resulted in a linear system for the source terms.  The right-hand side of this linear system is non-zero only in the presence of coalescence or non-linear evaporation.  The coefficient matrix depends on the choice of moments used in DQMOM.

We have compared this method, as well as the solution obtained with another Eulerian method:  the multi-fluid model, to the reference solution produced by a classical Lagrangian solver.
As far as coalescence phenomena are concerned, the efficiency of DQMOM has been shown to be better than the multi-fluid model due to its limited numerical diffusion in the size phase space, especially for the bimodal distribution function. However, as far as the evaporation process is concerned, it is
comparable to the multi-fluid model, but still needs a further study in order
to fully understand how to treat optimally the issue of the evaporative flux due to droplet disappearance. 
Although this issue has been so far 
neglected in the literature on moment methods, our study illustrates that it has an important effect on the moment dynamics.

The principal conclusion from this study is that DQMOM is numerically robust and straightforward to implement for the Williams spray equation and that it will be a very good candidate for more complex two-phase combustion
applications once the issue of the evaporative flux is further improved.

\section{Acknowledgments}

We would like to thank P.\ Villedieu for providing valuable help for the Lagrangian solver and several helpful discussions.
The present research was
done thanks to
a Young Investigator Award from the French Ministry of Research (New Interfaces of Mathematics - M.\ Massot, 2003-2006), the support of
the French Ministry of Research (Direction of the Technology) in the program:
``Recherche A\'eronautique sur le Supersonique'' (Project coordinator: M.\ Massot 2003-2006) and
an ANR (National Research Agency - 
France) Young Investigator Award (M.\ Massot, 2006-2009).
One of the author (R.\ O.\ Fox) was partially supported by the U.\ S.\ National Science Foundation (CTS-0403864).
The strong support of the
Ecole Centrale Paris through the invitation of R.\ O.\ Fox to Visiting Professor positions in 2005 and 2006 is gratefully acknowledged.

\bibliographystyle{plain}
\bibliography{references,biblio3}

\newpage

%

\end{document}